\documentclass{amsart}
\usepackage{amsfonts}
\usepackage{fancyhdr}
\usepackage[latin1]{inputenc}
\usepackage{amsmath}
\usepackage{thmdefs}
\usepackage{natbib}
\usepackage{graphicx}
\usepackage{hyperref}
\usepackage{amssymb}
\usepackage{dsfont}

\theoremstyle{definition} \theoremstyle{remark}
\numberwithin{equation}{section}

\renewcommand{\cite}{\citet}

\begin{document}

\author{Michel BRONIATOWSKI$^{*}$ and Amor KEZIOU$^{**}$}
\address{$^{*}$LSTA-Universit\'{e} Paris 6. e-mail: michel.broniatowski@upmc.fr
\newline
$^{**}$Laboratoire de Math\'ematiques (FRE 3111) CNRS,
Universit\'e de Reims  and LSTA-Universit\'e Paris 6. \\
 e-mail: amor.keziou@upmc.fr}
\date{October 2008}
\title[linear constraints with unknown parameters]{Estimation and tests for
models satisfying linear constraints with unknown parameter}
\maketitle

\begin{abstract}
We introduce estimation and test procedures through divergence
minimization for models satisfying linear constraints with unknown
parameter. Several statistical examples and motivations are given.
These procedures extend the empirical likelihood (EL) method and
share common features with generalized empirical likelihood (GEL).
We treat the problems of existence and characterization of the
divergence projections of probability measures on sets of signed
finite measures. Our approach allows for a study of the estimates
under misspecification. The asymptotic behavior of the proposed
estimates are studied using the dual representation of the
divergences and the explicit forms of the divergence projections.
We discuss the problem of the choice of the divergence under
various respects. Also we handle efficiency and robustness
properties of minimum divergence estimates. A simulation study
shows that the Hellinger divergence enjoys good efficiency and
robustness properties. \vspace{2mm} 
\end{abstract}




\section{Introduction and notation}

\noindent A model satisfying partly specified linear parametric
constraints is a family of distributions $\mathcal{M}^{1}$ all
defined on a same measurable space
$\left(\mathcal{X},\mathcal{B}\right)$, such that, for all $Q$ in
$\mathcal{M}^1$, the following condition holds
\begin{equation*}
\int g(x,\theta)~dQ(x)=0.
\end{equation*}
The unspecified parameter $\theta$ belongs to $\Theta$, an open
set  in $\mathbb{R}^{d}$. The function $g:={(g_{1},...,g_{l})}^T$
is defined on $\mathcal{X}\times \Theta$ with values in
$\mathbb{R}^{l}$, each of the $g_{i}$'s being real valued and the
functions $g_1,\ldots,g_l, \mathds{1}_\mathcal{X}$ are assumed
linearly independent. So $\mathcal{M}^{1}$ is defined through
$l$-linear constraints indexed by some $d-$dimensional parameter
$\theta$. Denote $M^{1}$ the collection of all probability
measures on $\left(\mathcal{X},\mathcal{B}\right)$, and
\begin{equation*}
\mathcal{M}_{\theta}^{1}:=\left\{ Q\in M^{1}~\text{ such that }
~\int g(x,\theta )~dQ(x)=0\right\}
\end{equation*}
so that
\begin{equation}
\mathcal{M}^{1}=\bigcup_{\theta \in \Theta }\mathcal{M}_{\theta
}^{1}. \label{modele}
\end{equation}
Assume now that we have at hand a sample $X_{1},...,X_{n}$ of
independent random variables (r.v.'s) with common unknown
distribution $P_0$. When $P_0$ belongs to the model
(\ref{modele}), we denote $\theta_0$ the value of the parameter
$\theta$ such that $\mathcal{M}_{\theta_0}$ contains $P_0$.
Obviously, we assume that $\theta_0$ is unique.\\

\noindent The scope of this paper is to propose new answers for
the classical following problems

\textit{Problem 1}: Does $P_0$ belong to the model
$\mathcal{M}^{1}$ ?

\textit{Problem 2}: When $P_0$ is in the model, which is the value
$\theta _{0}$ of the parameter for which $\int
 g(x,\theta_{0})~dP_0(x)=0$ ? Also can we perform simple and
composite tests for $\theta_0$ ? Can we construct confidence areas
for $\theta_{0}$ ? Can we give more efficient estimates for the
distribution function than the usual empirical
cumulative distribution function (c.d.f.) ?\\

\noindent We present some examples and motivations for the model
$(\ref{modele})$ and Problems 1 and 2.

\subsection{Statistical examples and motivations}

\begin{example}
Suppose that $P_0$ is the distribution of a pair of random
variables $(X,Y)$ on a product space
$\mathcal{X}\times\mathcal{Y}$ with known marginal distributions
$P_1$ and $P_2$. \cite{BickelRiovWellner1991} study efficient
estimation of $~\theta=\int h(x,y)~dP_0(x,y)$ for specified
function $h$. This problem can be handled in the present context
when the spaces $\mathcal{X}$ and $\mathcal{Y}$ are discrete and
finite. Denote $\mathcal{X}=\left\{x_1,\ldots,x_k\right\}$ and
$\mathcal{Y}=\left\{y_1,\ldots,y_r\right\}$. Consider an i.i.d.
bivariate sample $(X_i,Y_i), 1\leq i\leq n$ of the bivariate
random variable $(X,Y)$. The space $\mathcal{M}_\theta$ in this
case is the set of all p.m.'s  $Q$ on
$\mathcal{X}\times\mathcal{Y}$ satisfying $\int
g(x,y,\theta)~dQ(x,y)=0$ where
$g=(g_1^{(1)},\ldots,g_k^{(1)},g_1^{(2)},\ldots,g_r^{(2)},g_1)^T$,
$g_i^{(1)}(x,y,\theta)=\mathds{1}_{\left\{x_i\right\}\times
\mathcal{Y}}(x,y)-P_1(x_i)$,
$g_j^{(2)}(x,y,\theta)=\mathds{1}_{\mathcal{X}\times\left\{y_j\right\}}(x,y)
-P_2(y_j)$ for all $(i,j)$$\in$$\left\{1,\ldots,
k\right\}\times\left\{1,\ldots, r\right\}$, and
$g_1(x,y,\theta)=h(x,y)-\theta$. \textit{ Problem 1 } turns to be
the test for ``$P_0$ belongs to
$\bigcup_{\theta\in\Theta}\mathcal{M}_\theta$'', while \textit{
Problem 2} pertains to the estimation and tests for specific
values of $\theta$. Motivation and references for this problem are
given in \cite{BickelRiovWellner1991}.
\end{example}

\begin{example}(Generalized linear models). ~ Let $Y$ be a random variable and
$X$ a $l$-dimensional random vector. $Y$ and $X$ are linked
through
\begin{equation*}
Y=m(X,\theta_0)+\varepsilon
\end{equation*}
in which $m(.,.)$ is some specified real valued function and
$\theta_0$, the parameter of interest, belongs to some open set
$\Theta\subset\mathbb{R}^d$.~  $\varepsilon$ is a measurement
error. Denote $P_0$ the law of the vector variable $(X,Y)$ and
suppose that the true value $\theta_0$ satisfies the orthogonality
condition
\begin{equation*}
\int x(y-m(x,\theta_0))~dP_0(x,y)=0.
\end{equation*}
Consider an i.i.d. sample $(X_i,Y_i),~1\leq i\leq n$ of r.v.'s
with same distribution as $(X,Y)$. The existence of some
$\theta_0$ for which the above condition holds is given as the
solution of \textit{Problem 1}, while \textit{ Problem 2 } aims to
provide its explicit value; here $\mathcal{M}^1_\theta$ is  the
set of all p.m.'s $Q$ on $\mathbb{R}^{l+1}$ satisfying $\int
g(x,y,\theta)~dQ(x,y)=0$ with $g(x,y,\theta)=x(y-m(x,\theta))$.\\
\end{example}
\noindent \cite{Qin-Lawless1994} introduce various interesting
examples when (\ref{modele}) applies. In their example 1, they
consider the existence and estimation of the expectation $\theta$
of some r.v. $X$ when $E(X^2)=m(\theta)$ for some known function
$m(.)$. Another example is when a bivariate sample
$\left(X_i,Y_i\right)$ of i.i.d. r.v.'s is observed, the
expectation of $X_i$ is known and we intend to estimate $E(Y_i)$.
\cite{Haberman1984} and \cite{Sheehy1987} consider estimation of
$F(x)$ based on i.i.d. sample $X_1,\ldots,X_n$ with distribution
function $F$ when it is known that $\int T(x)~dF(x)=a$, for some
specified function $T(\cdot)$. For this problem, the function
$g(x,\theta)$ in the model (\ref{modele}) is equal to
$T(x)-\theta$ where $\theta=a$ is known. This example with $a$
unknown  is treated in details in Section 3 of the present paper.
We refer to \cite{Owen2001} for more statistical examples when
model (\ref{modele}) applies.\\

\noindent Another motivation for our work stems from confidence
region (C.R.) estimation techniques. The empirical likelihood
method provides such estimation (see \cite{Owen1990}). We will
extend this approach providing a wide range of such C.R.'s, each
one depending upon a specific criterion, one of those leading to
Owen's C.R.\\

\noindent An important estimator of $\theta_0$ is the generalized
method of moments (GMM) estimator of \cite{Hansen1982}. The
empirical likelihood approach  developed by \cite{Owen1988} and
\cite{Owen1990} has been adapted in the present setting by
\cite{Qin-Lawless1994} and \cite{Imbens1997} introducing the
empirical likelihood estimator (EL). The recent literature in
econometrics focusses on such models, the paper by
\cite{NeweySmith2003} provides an exhaustive list of works dealing
with the statistical properties of  GMM and generalized empirical
likelihood (GEL) estimators.\\

\noindent Our interest also lays in the behavior of the estimates
under misspecification. In the context of tests of hypothesis, the
statistics to be considered is some estimate of some divergence
between the unknown distributions of the data and the model. We
are also motivated by the behavior of those statistics under
misspecification, i.e., when the model is not appropriated to the
data. Such questions have not been addressed until now for those
problems in the general context of divergences.
\cite{Schennach2007} consider the asymptotic properties of the
empirical likelihood estimate under misspecification. As a by
product, we will prove that our proposal leads to consistent test
procedures; furthermore, the asymptotic behavior of the
statistics, under $\mathcal{H}_1$, provides the fundamental tool
in order to achieve Bahadur efficiency calculations (see
\cite{Nikitin1995}).\\

\noindent An important result due to \cite{NeweySmith2003} states
that EL estimate enjoys optimality properties in term of
efficiency when bias corrected among all GEL and GMM estimators.
Also \cite{Corcoran1998} and \cite{Baggerly1998} proved that in a
class of minimum discrepancy statistics, EL ratio is the only that
is Bartlett correctable. However, these results do not consider
the optimality properties of the tests for \textit{Problems} 1 and
2. Also, in connection with estimation problem, they do not
consider the properties of EL estimate with respect to robustness.
So, the question regarding divergence-based methods remains open
at least
in these two instances.\\

\noindent The approach which we develop is based on
\textit{minimum descrepancy estimates}, which have common features
with minimum distance techniques, using merely divergences. We
present  wide sets of estimates, simple and composite tests and
confidence regions for the parameter $\theta_{0}$ as well as
various test statistics for \textit{Problem 1}, all depending on
the choice of the divergence. Simulations show that the approach
based on Hellinger divergence enjoys good robustness and
efficiency properties when handling \textit{Problem 2}. As
presented in Section 5, empirical likelihood methods appear to be
a special case of the present approach.

\subsection{Minimum divergence estimates}

We first set some general definition and notation. Let $P$ be some
probability measure (p.m.). Denote $M^{1}(P)$ the subset of all
p.m.'s which are absolutely continuous (a.c.) with respect to $P$.
Denote $M$ the space of all signed finite measures on
$\left(\mathcal{X},\mathcal{B}\right)$ and $M(P)$ the subset  of
all signed finite measures a.c. w.r.t. $P$. Let $\varphi$ be a
convex function from $[-\infty, +\infty]$ onto $[0,+\infty]$ with
$\varphi(1)=0$. For any signed finite measure $Q$ in $M(P)$, the
$\phi-$divergence between $Q$ and the p.m. $P$ is defined through
\begin{equation}
\phi(Q,P):=\int \varphi\left(\frac{dQ}{dP}\right)~dP.
\label{divRusch}
\end{equation}
When $Q$ is not a.c. w.r.t. $P$, we set $\phi(Q,P)=+\infty$. This
definition extends \cite{Ruschendorf1984}'s one which applies for
$\phi-$divergences between p.m.'s; it also differs from
\cite{Csiszar1963}'s one, which requires a common dominating
$\sigma-$finite measure, noted $\lambda$, for $Q$ and $P$. Since
we will consider subsets of $M^{1}(P)$ and subsets of $M(P)$, it
is more adequate for our sake to use the definition
(\ref{divRusch}). Also note that all the just mentioned
definitions of $\phi-$divergences coincide on the set of all
p.m.'s a.c. w.r.t. $P$ and dominated by $\lambda$.\\

\noindent For all p.m. $P$, the mappings $Q\in
M\rightarrow\phi(Q,P)$ are convex and take nonnegative values.
When $Q=P$ then $\phi(Q,P)=0$. Further, if the function
$x\rightarrow\varphi(x)$ is strictly convex on neighborhood of
$x=1$, then the following basic property holds
\begin{equation}
\phi(Q,P)=0~\text{ if and only if }~Q=P.  \label{p.f.}
\end{equation}
\noindent All these properties are presented in
\cite{Csiszar1963}, \cite{Csiszar1967} and \cite{Liese-Vajda1987}
Chapter 1, for $\phi-$divergences defined on the set of all p.m.'s
$M^{1}$. When the $\phi$-divergences are defined on $M$, then the
same arguments as developed on $M^{1}$
hold.\\

\noindent When defined on $M^{1}$, the Kullback-Leibler $(KL)$,
modified Kullback-Leibler $(KL_{m})$, $\chi^{2}$, modified
$\chi^{2}$ $(\chi_{m}^{2})$, Hellinger $(H)$, and $L^{1}$
divergences are respectively associated to the convex functions
$\varphi(x)=x\log x-x+1$, $\varphi(x)=-\log x+x-1$,
$\varphi(x)=\frac{1}{2}{(x-1)}^{2}$,
$\varphi(x)=\frac{1}{2}{(x-1)}^{2}/x$,
$\varphi(x)=2{(\sqrt{x}-1)}^{2}$ and $\varphi(x)=\left\vert
x-1\right\vert$. All those divergences except the $L^{1}$ one,
belong to the class of power divergences introduced in
\cite{Cressie-Read1984} (see also \cite{Liese-Vajda1987} Chapter
2). They are defined through the class of convex functions
\begin{equation}\label{gamma convex functions}
x\in \mathbb{R}_{+}^*
\mapsto\varphi_{\gamma}(x):=\frac{x^{\gamma}-\gamma
x+\gamma-1}{\gamma(\gamma-1)}
\end{equation}
if $\gamma\in\mathbb{R}\setminus \left\{0,1\right\}$ and by
$\varphi_{0}(x):=-\log x+x-1$ and $\varphi_{1}(x):=x\log x-x+1$.
For all $\gamma\in\mathbb{R}$,
$\varphi_\gamma(0):=\lim_{x\downarrow 0}\varphi_\gamma (x)$ and
$\varphi_\gamma(+\infty):=\lim_{x\uparrow +\infty}\varphi_\gamma
(x)$. So, the $KL-$divergence is associated to $\varphi_1$, the
$KL_m$ to $\varphi_0$, the $\chi^2$ to $\varphi_2$, the $\chi^2_m$
to $\varphi_{-1}$ and the Hellinger distance to $\varphi_{1/2}$.
For all $\gamma\in\mathbb{R}$, we sometimes denote $\phi_\gamma$
the divergence associated to the convex function $\varphi_\gamma$.
We define the derivative of $\varphi_\gamma$ at $0$ by
$\varphi_\gamma'(0):=\lim_{x\downarrow 0}\varphi_\gamma'(x)$. We
extend the definition of the power divergences functions $Q\in
M^{1}\rightarrow \phi_{\gamma}(Q,P)$ onto the whole set of signed
finite measures $M$ as follows. When the function $x\rightarrow
\varphi_{\gamma }(x)$ is not defined on $(-\infty, 0[$ or when
$\varphi_{\gamma}$ is defined on $\mathbb{R}$ but is not a convex
function we extend the definition of $\varphi_{\gamma}$ through
\begin{equation} \label{gamma convex functions sur R}
x\in
[-\infty,+\infty]\mapsto\varphi_\gamma(x)\mathds{1}_{[0,+\infty]}(x)+
(\varphi_\gamma'(0)x+\varphi_\gamma(0))\mathds{1}_{[-\infty,0[}(x).
\end{equation}
For any convex function $\varphi$, define the \textit{domain} of
$\varphi$ through
\begin{equation}
D_\varphi=\left\{ x\in [-\infty,+\infty]~\text{ such that
}~\varphi (x)<+\infty \right\}. \label{domaine de varphi}
\end{equation}
Since $\varphi$ is convex, $D_\varphi$ is an interval which may be
open or not, bounded or unbounded. Hence, write $D_\varphi:=(a,b)$
in which $a$ and $b$ may be finite or infinite. In this paper, we
will only consider $\varphi$ functions defined on
$[-\infty,+\infty]$ with values in $[0, +\infty]$ such that
$a<1<b$, and which satisfy $\varphi(1)=0$, are strictly convex and
are $\mathcal{C}^2$ on the interior of its domain $D_\varphi$; we
define $\varphi(a)$, $\varphi'(a)$, $\varphi''(a)$, $\varphi(b)$,
$\varphi'(b)$ and $\varphi''(b)$ respectively by
$\varphi(a):=\lim_{x \downarrow a}\varphi(x)$,
$\varphi'(a):=\lim_{x \downarrow a}\varphi'(x)$,
$\varphi''(a):=\lim_{x \downarrow a}\varphi''(x)$,
$\varphi(b):=\lim_{x \uparrow b}\varphi(x)$,
 $\varphi'(b):=\lim_{x \uparrow b}\varphi'(x)$
and $\varphi''(b):=\lim_{x \uparrow b}\varphi''(x)$. These
quantities may be finite or infinite. All the functions
$\varphi_\gamma$ (see (\ref{gamma convex functions sur R}))
satisfy these conditions.

\begin{definition}
Let $\Omega$ be some subset in $M$. The $\phi-$divergence between
the set $\Omega$ and a p.m. $P$, noted $\phi(\Omega,P)$,  is
\begin{equation*}
\phi(\Omega,P):=\inf_{Q\in\Omega}\phi(Q,P).
\end{equation*}
\end{definition}

\begin{definition}
Assume that $\phi(\Omega,P)$ is finite. A measure $Q^*\in\Omega$
such that
\begin{equation*}
    \phi(Q^*,P)\leq \phi(Q,P) ~\text{ for all }~Q\in\Omega
\end{equation*}
is called a $\phi-$projection of $P$ onto $\Omega$. This
projection may not exist, or may be not defined uniquely.
\end{definition}
\noindent We will make use of the concept of $\phi-$divergences in
order to perform estimation and tests for the model
(\ref{modele}). So, let $X_{1},...,X_{n}$ denote an i.i.d. sample
of r.v.'s with common distribution $P_{0}$. Let $P_{n}$ be the
empirical measure pertaining to this sample, namely
\begin{equation*}
P_{n}:=\frac{1}{n}\sum_{i=1}^{n}\delta_{X_{i}}
\end{equation*}
in which $\delta_{x}$ is the Dirac measure at point $x$. When
$P_{0}$ and all $Q\in\mathcal{M}^1$ share the same discrete finite
support $S$, then the $\phi$-divergence $\phi(Q,P_{0})$ can be
written as
\begin{equation}\label{div discrete}
\phi(Q,P_{0})=\sum_{j\in S}\varphi
\left(\frac{Q(j)}{P_{0}(j)}\right)P_{0}(j).
\end{equation}
In this case, $\phi(Q,P_{0})$ can be estimated simply through the
plug-in of $P_{n}$ in (\ref{div discrete}), as follows
\begin{equation}\label{div est support fini}
\widehat{\phi }(Q,P_{0}):=\sum_{j\in S}\varphi \left(
\frac{Q(j)}{P_{n}(j)}\right) P_{n}(j).
\end{equation}
In the same way, for any $\theta$ in $\Theta$,
$\phi\left(\mathcal{M}_{\theta}^{1},P_{0}\right)$ is estimated by
\begin{equation}\label{estim de phi M1}
\widehat{\phi}\left(\mathcal{M}_{\theta }^{1},P_{0}\right)
:=\inf_{Q\in \mathcal{M}_{\theta}^{1}}\sum_{j\in S}\varphi \left(
\frac{Q(j)}{P_{n}(j)}\right) P_{n}(j),
\end{equation}
and $\phi\left(\mathcal{M}^1,P_0\right)=\inf_{\theta\in\Theta}
\phi\left(\mathcal{M}_\theta^1,P_0\right)$ can be estimated by
\begin{equation}\label{estim de inf de phi M1}
    \widehat{\phi}\left(\mathcal{M}^1,P_0\right):=\inf_{\theta\in\Theta}
    \inf_{Q\in \mathcal{M}_{\theta}^{1}}\sum_{j\in S}\varphi \left(
\frac{Q(j)}{P_{n}(j)}\right) P_{n}(j).
\end{equation}
 By uniqueness of $\inf_{\theta \in \Theta }\phi \left(
\mathcal{M}_{\theta}^{1},P_{0}\right) $ and since this infimum is
reached at $\theta =\theta_{0}$, we estimate $\theta_{0}$ through
\begin{equation}\label{estim de theta phi M1}
\widehat{\theta}_\phi:=\arg \inf_{\theta\in\Theta} \inf_{Q\in
\mathcal{M}_{\theta}^{1}}\sum_{j\in S}\varphi \left(
\frac{Q(j)}{P_{n}(j)}\right) P_{n}(j).
\end{equation}
The infimum in (\ref{estim de phi M1}) (i.e., the projection of
$P_n$ on $\mathcal{M}_\theta^1$) may be achieved on the frontier
of $\mathcal{M}_{\theta }^{1}$. In this case the Lagrange method
is not valid. Hence, we endow our statistical approach in the
global context of signed finite measures with total mass 1
satisfying the linear constraints.
\begin{equation}\label{modele signees finies}
\mathcal{M}_{\theta }:=\left\{ Q\in M~\text{ such that }~\int
dQ=1\text{ and }\int g(x,\theta )~dQ(x)=0\right\}
\end{equation}
and
\begin{equation}
\mathcal{M}:=\bigcup_{\theta\in\Theta}\mathcal{M}_{\theta},
\end{equation}
sets of signed finite measures that replace $\mathcal{M}_{\theta
}^{1}$ and $\mathcal{M}^1$.

As above, we estimate $\phi(\mathcal{M}_\theta,P_0)$,
$\phi(\mathcal{M},P_0)$ and $\theta_0$ respectively by
\begin{equation}\label{estim de phiMtheta cas disceret}
 \widehat{\phi}\left( \mathcal{M}_{\theta
},P_{0}\right):=\inf_{Q\in \mathcal{M}_{\theta}}\sum_{j\in
S}\varphi \left( \frac{Q(j)}{P_{n}(j)}\right) P_{n}(j),
\end{equation}
\begin{equation}\label{estim de phiM cas discret}
 \widehat{\phi}\left(\mathcal{M},P_0\right):=\inf_{\theta\in\Theta}
 \inf_{Q\in \mathcal{M}_{\theta}}\sum_{j\in
 S}\varphi \left( \frac{Q(j)}{P_{n}(j)}\right) P_{n}(j),
\end{equation}
and
\begin{equation}\label{estim de theta0 cas discret}
 \widehat{\theta}_\phi:=\arg \inf_{\theta \in \Theta } \inf_{Q\in \mathcal{M}_{\theta}}
 \sum_{j\in S}\varphi \left( \frac{Q(j)}{P_{n}(j)}\right)
 P_{n}(j).
\end{equation}
Enhancing $\mathcal{M}^1$ to $\mathcal{M}$ is motivated by the
following arguments
\begin{enumerate}
 \item [-] For all $\theta$ in $\Theta$, denote $Q_1^*$ and $Q^*$  respectively the
 projection of $P_n$ on $\mathcal{M}_\theta^1$ and on
 $\mathcal{M}_\theta$,  as defined in
 (\ref{estim de phi M1}) and in (\ref{estim de phiMtheta cas disceret}).
 If $Q_1^*$ is an interior point of
 $\mathcal{M}_\theta^1$, then, by Proposition
 \ref{caracterisation de Q star n} below, it coincides with $Q^*$,
 the projection of $P_n$ on
 $\mathcal{M}_\theta$, i.e., $Q_1^*=Q^*$. Therefore, in this case, both approaches
 coincide.
 \item [-] It may occur that for some $\theta$ in $\Theta$, $Q_1^*$, the
 projection  of $P_n$ on $\mathcal{M}_\theta^1$, is a frontier
 point of $\mathcal{M}_\theta^1$, which makes a real difficulty
 for the estimation procedure. We will prove in Theorem \ref{loi limite de thetan et cn}
 that $\widehat{\theta}_\phi$, defined in (\ref{estim de theta0 cas discret})
 and which replaces (\ref{estim de theta phi M1}), converges to
 $\theta_0$. This validates the substitution of the sets $\mathcal{M}_\theta^1$
 by the sets $\mathcal{M}_\theta$. In the context of a test problem,
 we will prove that the asymptotic distributions of the test
 statistics pertaining to Problem 1 and 2 are unaffected by this
 change.
\end{enumerate}
This modification motivates the above extensions in the
 definitions of the $\varphi$ functions on $[-\infty,+\infty]$ and of the
 $\phi$-divergences on the whole space of finite signed measures $M$.

\vskip 0.5cm

\noindent In the case when $Q$ and $P_0$ share different discrete
finite support or share  same or different discrete infinite or
continuous support, then formula (\ref{div est support fini}) is
not defined, due to lack of absolute continuity of $Q$ with
respect to $P_{n}$. Indeed
\begin{equation}\label{perte c.a.}
\widehat{\phi}(Q,P_{0}):=\phi(Q,P_{n})=+\infty.
\end{equation}
The plug-in estimate of $\phi(\mathcal{M}_\theta,P_0)$ is
\begin{equation}\label{plug-in est of phi}
\widehat{\phi}(\mathcal{M}_\theta,P_0):=\inf_{Q\in\mathcal{M}_\theta}\phi(Q,P_n)=
\inf_{Q\in\mathcal{M}_\theta}\int
\varphi\left(\frac{dQ}{dP_n}(x)\right)~dP_n(x).
\end{equation}
If  the infimum exists, then it is clear that it is reached at a
signed finite measure (or probability measure) which is  a.c.
w.r.t. $P_n$. So, define the sets
\begin{equation}\label{m theta n}
\mathcal{M}_\theta^{(n)}:=\left\{Q\in M ~\text{ such that }~Q\ll
P_n, ~ \sum_{i=1}^{n}Q(X_i)=1 ~\text{ and
}~\sum_{i=1}^{n}Q(X_i)g(X_i,\theta)=0 \right\},
\end{equation}
which may be seen as subsets of $\mathbb{R}^n$. Then, the plug-in
estimate (\ref{plug-in est of phi}) of
$\phi(\mathcal{M}_\theta,P_0)$ can be written as
\begin{equation}   \label{estim de phiMtheta cas cont}
   \widehat{\phi}(\mathcal{M}_\theta,P_0)=\inf_{Q\in\mathcal{M}_\theta^{(n)}}
   \frac{1}{n}\sum_{i=1}^{n}\varphi\left(nQ(X_i)\right).
\end{equation}
In the same way,
$\phi(\mathcal{M},P_0):=\inf_{\theta\in\Theta}\inf_{Q\in\mathcal{M}_\theta}
\phi(Q,P_0)$ can be estimated by
\begin{equation}\label{estim de phiM  cas cont}
  \widehat{\phi}(\mathcal{M},P_0)=\inf_{\theta\in\Theta}
  \inf_{Q\in\mathcal{M}_\theta^{(n)}}\frac{1}{n}\sum_{i=1}^{n}\varphi\left(nQ(X_i)\right).
\end{equation}
By uniqueness of
$\inf_{\theta\in\Theta}\phi(\mathcal{M}_\theta,P_0)$ and since
this infimum is reached at $\theta=\theta_0$, we estimate
$\theta_0$ through
\begin{equation}\label{estim de theta0 cas cont}
    \widehat{\theta}_\phi=\arg\inf_{\theta\in\Theta}\inf_{Q\in\mathcal{M}_\theta^{(n)}}
    \frac{1}{n}\sum_{i=1}^{n}\varphi\left(nQ(X_i)\right).
\end{equation}
Note that, when $P_0$ and all $Q\in\mathcal{M}^1$ share the same
discrete finite support, then the estimates (\ref{estim de theta0
cas cont}), (\ref{estim de phiM cas cont}) and (\ref{estim de
phiMtheta  cas cont}) coincide respectively with (\ref{estim de
theta0 cas discret}), (\ref{estim de phiM cas discret}) and
(\ref{estim de phiMtheta cas disceret}). Hence, in the sequel, we
study the estimates $\widehat{\phi}(\mathcal{M}_\theta, P_0)$,
$\widehat{\phi}(\mathcal{M}, P_0)$ and $\widehat{\theta}_\phi$ as
defined in (\ref{estim de phiMtheta  cas cont}), (\ref{estim de
phiM  cas cont})    and (\ref{estim de theta0 cas cont}),
respectively. We propose to call the estimates
$\widehat{\theta}_\phi$ defined in (\ref{estim de theta0 cas
cont}) ``Minimum Empirical $\phi$-Divergences Estimates''
(ME$\phi$DE's). As will be noticed later on, the empirical
likelihood paradigm (see \cite{Owen1988} and \cite{Owen1990}),
which is based on this plug-in approach, enters as a special case
of the statistical issues related to estimation and tests based on
$\phi-$divergences with $\varphi(x)=\varphi_{0}(x)=-\log x+x-1$,
namely on $KL_m-$divergence. The empirical log-likelihood ratio
for the model (\ref{modele signees finies}), in the context of
$\phi$-divergences, can be written as
$-n\widehat{KL_m}(\mathcal{M}_\theta,P_0)$. In the case of a
single functional, for example when $g(x,\theta)=x-\theta$ with
$x$ and $\theta$ belong to $\mathbb{R}$, \cite{Owen1988} shows
that $2n\widehat{KL_m}(\mathcal{M}_\theta,P_0)$ has an asymptotic
$\chi^2_{(1)}$ distribution when $P_0$ belongs to
$\mathcal{M}_\theta$. (see \cite{Owen1988} Theorem 1). This result
is a nonparametric version of Wilks's theorem (see
\cite{Wilks1938}). In the multivariate case,  the same result
holds (see \cite{Owen1990} Theorem 1). When we want to extend the
arguments used in \cite{Owen1988} and \cite{Owen1990}  in order to
study the limiting behavior of the statistics
$\widehat{\phi}(\mathcal{M}_\theta,P_0)$,  when
$P_0\not\in\mathcal{M}_\theta$ (for example, when
$\theta_0\neq\theta$), most limiting arguments become untractable.
We propose to use the so-called ``dual representation of
$\phi-$divergences''(see \cite{Keziou2003}), a device which is
well known for the Kullback-Leibler divergence in the context of
large deviations, and which has been used in parametric statistics
in \cite{Keziou2003} and \cite{Broniatowski-Keziou2008}. The
estimates then turn to be M-estimates whose limiting distributions
are obtained through classical methods. On the other hand, the
obtention of the limit distributions of the statistics
$\widehat{\phi}(\mathcal{M}_\theta,P_0)$ when $P_0\not\in
\mathcal{M}_\theta$, requires the study of the existence and the
characterization of the projection of the p.m. $P_0$ on the sets
$\mathcal{M}_\theta$.\\

\noindent This paper is organized as follows : In Section 3,  we
study the asymptotic behavior of the proposed estimates
(\ref{estim de phiMtheta  cas cont}), (\ref{estim de phiM  cas
cont}) and (\ref{estim de theta0 cas cont}) giving solutions to
$\textit{Problem 2}$. We then address \textit{Problem 1}, namely :
does there exist some $\theta_{0}$ in $\Theta$ for which $P_0$
belongs to $\mathcal{M}_{\theta_{0}}$? In Section 4, extending the
result by \cite{Qin-Lawless1994}, we give new estimates for the
distribution function using the $\phi$-projections of $P_n$ on the
model $\mathcal{M}$. We show that the new estimates of the
distribution function are generally more efficient than the
empirical cumulative distribution function. Section 5 illustrates
the concept of empirical likelihood in the context of
$\phi$-divergences techniques. In Section 6, we focus on
robustness and efficiency of the ME$\phi$D estimates. A simulation
study aims at emphasizing the specific advantage of the choice of
the Hellinger divergence in relation with robustness and
efficiency considerations. All proofs are in Section 7.

\section{Estimation for Models  satisfying
 Linear\\ Constraints}

\noindent At this point, we must introduce some notational
convention for sake of brevity and clearness. For any p.m. $P$ on
$\mathcal{X}$ and any measurable real function $f$ on
$\mathcal{X}$, $Pf$ denotes $\int f(x)~dP(x)$. For example,
$P_0g_j(\theta)$ will be used instead of $\int
g_j(\theta,x)~dP_0(x)$. Hence, we are led to define the following
functions : denote $\overline{g}$ the function defined on
$\mathcal{X}\times\Theta$ with values in $\mathbb{R}^{l+1}$ by
\begin{equation*}
\begin{array}{ccccl}
\overline{g} & : & \mathcal{X}\times \Theta   & \rightarrow &
\mathbb{R}^{(l+1)}
\\
 &  & (x,\theta) & \mapsto & \overline{g}(x,\theta):=
{\left(\mathds{1}_\mathcal{X}(x),g_1(x,\theta),\ldots,g_l(x,\theta)\right)}^T,
\end{array}
\end{equation*}
and for all $\theta\in\Theta$, denote also $\overline{g}(\theta)$,
$g(\theta)$, $g_j(\theta)$ the functions defined respectively by
\begin{equation*}
\begin{array}{ccccl}
\overline{g}(\theta) & : & \mathcal{X} & \rightarrow &
\mathbb{R}^{l+1}\\
& & x & \mapsto &
\overline{g}(x,\theta):={\left(g_0(x,\theta),g_1(x,\theta),
\ldots,g_l(x,\theta)\right)}^T,\text{ where
}g_0(x,\theta):=\mathds{1}_\mathcal{X}(x),\\
g(\theta) & : & \mathcal{X} & \rightarrow &
\mathbb{R}^l\\
& & x & \mapsto &
g(x,\theta):={\left(g_1(x,\theta),\ldots,g_l(x,\theta)\right)}^T\\
\text{and} & & & &\\
g_j(\theta) & : & \mathcal{X} & \rightarrow &
\mathbb{R}\\
& & x & \mapsto & g_j(x,\theta),~~\text{ for all }
j\in\left\{0,1,\ldots,l\right\}.\\\\
\end{array}
\end{equation*}

\noindent We now turn back to the setting defined in the
Introduction and consider model (\ref{modele signees finies}). For
fixed $\theta$ in $\Theta$, define the class of functions
\begin{equation*}
\mathcal{F}_{\theta
}:=\left\{g_{0}(\theta),g_{1}(\theta),\ldots,g_{l}(\theta)\right\},
\end{equation*}
and consider the set of finite signed measures
$\mathcal{M}_{\theta}$ defined by $(l+1)$ linear constraints as
defined in $\left(\ref{modele signees finies}\right)$
\begin{equation*}
\mathcal{M}_{\theta }:=\left\{ Q\in
M_{\mathcal{F}_{\theta}}~\text{ such that }~\int dQ(x)=1 \text{
and }\int g(x,\theta )~dQ(x)=0\right\}.
\end{equation*}

\noindent We present explicit tractable conditions for the
estimates (\ref{estim de phiMtheta cas cont}), (\ref{estim de phiM
cas cont}) and (\ref{estim de theta0 cas cont}) to be well
defined. This will be done in  Propositions \ref{existence de Q
star n}, Remark \ref{condition d'existence de Qn cq}, Proposition
\ref{existence de Q star} and Remark
\ref{condition d'existence de Qstar cq} below.\\
First, we present sufficient conditions which assess the existence
of the infimum in (\ref{estim de phiMtheta cas cont}), noted
$\widehat{Q_\theta^*}$, the projection of $P_n$ on
$\mathcal{M}_\theta$. We also provide conditions under which the
Lagrange method can be used to characterize
$\widehat{Q_\theta^*}$. The Fenchel-Legendre transform of
$\varphi$ will be denoted $\varphi ^{\ast }$, i.e.,
\begin{equation}
t\in \mathbb{R}\mapsto \varphi ^{\ast }(t):=\sup_{x\in
\mathbb{R}}\left\{ tx-\varphi (x)\right\}.
\end{equation}
Define
\begin{equation}\label{D phi n}
    \mathcal{D}_\phi^{(n)}:=\left\{Q\in M ~\text{ such that }~Q\ll P_n~
    \text{ and }
    \frac{1}{n}\sum_{i=1}^{n}\varphi\left(nQ(X_i)\right)<\infty~\right\},
\end{equation}
i.e., the \textit{domain} of the function
\begin{equation*}
\left(Q(X_1),\ldots,Q(X_n)\right)^T\in\mathbb{R}^n\mapsto
\frac{1}{n}\sum_{i=1}^{n}\varphi\left(nQ(X_i)\right).
\end{equation*}
We have
\begin{proposition}\label{existence de Q star n}
Assume that there exists some measure $R$ in the interior of
$\mathcal{D}_\phi^{(n)}$ and in $\mathcal{M}_\theta^{(n)}$ such
that for all $Q$ in $\partial \mathcal{D}^{(n)}_\phi$,  the
frontier of $\mathcal{D}^{(n)}_\phi$,
 we have
\begin{equation}\label{condition d'existence de Q star n}
   \frac{1}{n}\sum_{i=1}^{n}\varphi\left(nR(X_i)\right)<
   \frac{1}{n}\sum_{i=1}^{n}\varphi\left(nQ(X_i)\right).
\end{equation}
Then the following holds
\begin{enumerate}
\item[(i)] there exists an unique $\widehat{Q_\theta^*}$ in
$\mathcal{M}_\theta^{(n)}$ such that
     \begin{equation}\label{Q satr n est atteint}
          \inf_{Q\in\mathcal{M}_\theta^{(n)}}\frac{1}{n}\sum_{i=1}^{n}
          \varphi\left(nQ(X_i)\right)=\frac{1}{n}\sum_{i=1}^{n}
          \varphi\left(n\widehat{Q_\theta^*}(X_i)\right)
     \end{equation}
\item[(ii)] $\widehat{Q_\theta^*}$ is an interior point of
$\mathcal{D}_\phi^{(n)}$ and satisfies for all $i=1,\ldots,n$
\begin{equation}\label{caracterisation de Q star n}
    \widehat{Q_\theta^*}(X_i)=\frac{1}{n}\overleftarrow{\varphi'}
    \left(\sum_{j=0}^l\widehat{c}_jg_j(X_i,\theta)\right),
\end{equation}
where
$\left(\widehat{c}_0,\widehat{c}_1,\ldots,\widehat{c}_l\right)^T:=\widehat{c_\theta}$
is solution of the system of equations
\begin{equation} \label{syst empirique}
\left\{
\begin{array}{ccl}
\int \overleftarrow{\varphi ^{\prime }}\left(
\widehat{c}_{0}+\sum_{i=1}^{l}\widehat{c}_{i}g_{i}(x,\theta )\right) ~dP_n(x) & = & 1 \\
\int g_{j}(x,\theta )\overleftarrow{\varphi ^{\prime }}\left(
\widehat{c}_{0}+\sum_{i=1}^{l}\widehat{c}_{i}g_{i}(x,\theta)\right)
~dP_n(x) & = & 0,~~ j=1,\ldots,l.
\end{array}
\right.
\end{equation}
\end{enumerate}
\end{proposition}

\begin{example}
For the $\chi^{2}-$divergence, we have
$\mathcal{D}^{(n)}_{\chi^2}=\mathbb{R}^{n}$. Hence condition
(\ref{condition d'existence de Q star n}) holds whenever
$\mathcal{M}_\theta^{(n)}$ is not void. Therefore, the above
Proposition holds always independently upon the distribution
$P_0$. More generally, the above Proposition holds for any
$\phi$-divergence which is associated to $\varphi$ function
satisfying $D_\varphi=\mathbb{R}$. (See (\ref{domaine de varphi})
for the definition of $D_\varphi$).
\end{example}

\begin{example}
In the case of the modified Kullback-Leibler divergence, which
turns to coincide with the empirical likelihood technique (see
Section 5), we have $\mathcal{D}^{(n)}_{KL_m}= \left( ]0,+\infty[
\right) ^{n}$. For $\alpha$ in $\Theta$, define the assertion
\begin{equation}\label{cond owen}
\begin{array}{rcl}
\text{ there exists }~q & = & (q_{1},...,q_{n})\text{ in
}\mathbb{R}^{n}~\text{
with }~ 0<q_{i}<1\text{ for all }i=1,...,n   \\
\text{ and }~\sum_{i=1}^{n}q_{i}g_{j}(X_{i},\alpha ) & = &
0~\text{ for all } ~j=1,...,l.
\end{array}
\end{equation}
A sufficient condition, in order to assess that condition
(\ref{condition d'existence de Q star n}) in the above Proposition
holds, is when $\left( \ref{cond owen}\right)$ holds for
$\alpha=\theta$. In the case when $g(x,\theta )=x-\theta$, this is
precisely what is checked in (\cite{Owen1990}), p. 100, when
$\theta$ is an interior point of the convex hull of
$(X_{1},...,X_{n})$.
\end{example}

\begin{example}
For the modified $\chi^2-$divergence, we have
$\mathcal{D}^{(n)}_{\chi^2_m}=\left( ]0,\infty [ \right) ^{n}$,
and therefore, condition (\ref{cond owen}) for $\alpha=\theta$ is
sufficient for the condition (\ref{condition d'existence de Q star
n}) to holds. So, conditions which assess the existence of the
projection $\widehat{Q_\theta^*}$ are the same for the modified
$\chi^2-$divergence and the
$KL_m$-divergence.\\
\end{example}

\begin{remark}\label{condition d'existence de Qn cq}
If there exists some $Q_0\in\mathcal{M}_\theta^{(n)}$ such that
\begin{equation}\label{constraint qualification Pn}
a < \inf_i nQ_0(X_i)\leq \sup_i nQ_0(X_i) < b,
\end{equation}
then applying Corollary 2.6 in \cite{BorweinLewis1991}, we get
\begin{equation*}
 \inf_{Q\in\mathcal{M}_\theta^{(n)}}\frac{1}{n}\sum_{i=1}^{n}
          \varphi\left(nQ(X_i)\right)=\sup_{t\in \mathbb{R}^{(l+1)}}
          \left\{t_0-\int \psi\left(t^T\overline{g}(x,\theta)\right)~dP_n(x)\right\}
\end{equation*}
with dual attainement. Furthermore, if
\begin{equation*}
 \varphi'(a)< \inf_i {\widehat{c}_\theta}^T\overline{g}(X_i,\theta)\leq
 \sup_i {\widehat{c}_\theta}^T\overline{g}(X_i,\theta)<\varphi'(b),
\end{equation*}
with $\widehat{c}_\theta$ a dual optimal, then the unique
projection $\widehat{Q_\theta^*}$ of $P_n$ on
$\mathcal{M}_\theta^{(n)}$ is
given by (\ref{caracterisation de Q star n}).\\
\end{remark}
\noindent We will make use of the dual representation of
$\phi$-divergences (see \cite{Bronia_Kez2006_STUDIA} theorem 4.4).
So, define
\begin{equation}\label{C theta}
\mathcal{C}_{\theta }:=\left\{ t\in \mathbb{R}^{l+1}~\text{ such
that }~ t^T\overline{g}(.,\theta)~\text{ belongs to }~\text{ Im
}\varphi'~~(P_0-a.s.)\right\},
\end{equation}
and
\begin{equation}\label{C theta n}
\mathcal{C}_{\theta}^{(n)}:=\left\{ t\in\mathbb{R}^{l+1}~\text{
such }~ t^T\overline{g}(X_i,\theta)~\text{ belongs to }~\text{ Im
}\varphi' ~\text{ for all }~ i=1,\ldots,n\right\}.
\end{equation}
We omit the subscript $\theta$ when unnecessary. Note that both
$\mathcal{C}_\theta$ and $\mathcal{C}_\theta^{(n)}$ depend upon
the function
$\varphi$ but, for  simplicity,  we omit the subscript $\varphi$.\\

\noindent  If $P_0$ admits a  projection $Q_\theta^*$ on
$\mathcal{M}_{\theta }$ with the same support as $P_0$,
 using the second  part in Corollary 3.5 in \cite{Bronia_Kez2006_STUDIA}, there exist
constants $c_{0},\ldots,c_{l}$, obviously depending on $\theta$,
such that
\begin{equation*}
\varphi'\left(\frac{dQ^*_\theta}{dP_0}(x)\right)=c_{0}+\sum_{j=1}^{l}c_{j}g_{j}(x,\theta),
~ \text{ for all } ~x ~(P_0-a.s.).
\end{equation*}
Since $Q^*_\theta$ belongs to $\mathcal{M}_\theta$, the real
numbers $c_0,c_{1},\ldots ,c_{l}$  are solutions of
\begin{equation}\label{syst teta fixe}
\left\{
\begin{array}{ccl}
\int {\varphi^{\prime}}^{-1}\left(
c_{0}+\sum_{j=1}^{l}c_{j}g_{j}(x,\theta)\right)~dP_0(x) & = & 1 \\
\int g_{j}(x,\theta ){\varphi^{\prime}}^{-1}\left(
c_{0}+\sum_{j=1}^{l}c_{j}g_{j}(x,\theta)\right)~dP_0(x) & = & 0,~
j=1,\ldots ,l.
\end{array}
\right.
\end{equation}
Since $Q\mapsto \phi(Q,P_0)$ is strictly convex,  the projection
$Q_\theta^*$ of $P_0$ on the convex set $\mathcal{M}_\theta$ is
unique. This implies, by \cite{Bronia_Kez2006_STUDIA} Corollary
part 1, that the solution
\begin{equation*}
c_\theta:=(c_0,c_1,\ldots,c_l)^T
\end{equation*}
of the system (\ref{syst teta fixe}) is unique provided that the
functions $g_i(\theta)$ are linearly independent. Further, using
the dual representation of $\phi$-divergences (see
\cite{Bronia_Kez2006_STUDIA} Theorem 4.4), we get
\begin{equation*}
\phi(\mathcal{M}_{\theta},P_0):=\phi(Q_\theta^*,P_0)=\sup_{f\in
\mathcal{F}}\left\{\int f~dQ_\theta^*-\int
\varphi^*(f)~dP_0\right\},
\end{equation*}
and the sup is unique and is reached at
$f=\varphi^{\prime}(dQ_\theta^*/dP_0)=c_{0}+\sum_{j=1}^{l}c_{j}g_{j}(.,\theta)$,
if it belongs to $\mathcal{F}$. This motivates the choice of the
class $\mathcal{F}$ through
\begin{equation*}
\mathcal{F}:=\left\{ x\rightarrow t^T\overline{g}(x,\theta
)~~\text{ for }~~ t~~ \text{ in }~~\mathcal{C}_{\theta }\right\}.
\end{equation*}
It is the smallest class of functions that contains
$\varphi'(dQ_\theta^*/dP_0)$ and which does not presume any
knowledge on $Q_\theta^*$. We thus obtain
\begin{equation*}
\phi(\mathcal{M}_{\theta},P_0)=\sup_{t\in\mathcal{C}_{\theta}}\int
m(x,\theta,t)~dP_0(x),
\end{equation*}
where $m(\theta,t)$ is the function defined on $\mathcal{X}$ by
\begin{equation}\label{m x theta t}
\begin{array}{ccl}
 x\in\mathcal{X} & \mapsto & m(x,\theta,t):=
 t_0-\varphi^*\left(t^T\overline{g}(x,\theta)\right)=\nonumber\\
& & t_0-\left(t^T\overline{g}(x,\theta)\right) {\varphi
^{\prime}}^{-1}\left(t^T\overline{g}(x,\theta)\right) +\varphi
\left({\varphi^{\prime}}^{-1}\left(
t^T\overline{g}(x,\theta)\right)\right).
\end{array}
\end{equation}
With the above notation, we state
\begin{equation}\label{sup P0 m}
\phi(\mathcal{M}_{\theta},P_0)=\sup_{t\in\mathcal{C}_{\theta}}
P_0m(\theta,t).
\end{equation}
So, a natural estimate of $\phi(\mathcal{M}_\theta,P_0)$ is
\begin{equation}
\sup_{t\in\mathcal{C}_{\theta}^{(n)}} P_nm(\theta,t)
\end{equation}
which coincides with the estimate defined in (\ref{estim de
phiMtheta cas cont}). Hence, we can write
\begin{equation}\label{estim de
phiMtheta dual} \widehat{\phi}(\mathcal{M}_\theta,P_0)=\sup_{t\in
\mathcal{C}_\theta^{(n)}}P_nm(\theta,t).
\end{equation}
which transforms the constrained optimization in (\ref{estim de
phiMtheta cas cont}) into the above unconstrained one.\\

\noindent On the other hand, the sup in (\ref{sup P0 m})  is
reached at $t_0=c_{0},\ldots ,t_l=c_{l}$ which are solutions of
the system of equations $\left( \ref{syst teta fixe}\right)$,
i.e.,
\begin{equation}\label{c}
    c_\theta=\arg\sup_{t\in\mathcal{C}_\theta} P_0m(\theta,t).
\end{equation}
So, a natural estimate of $~c_\theta~$ in (\ref{c}) is therefore
defined through
\begin{equation}
\arg\sup_{t\in\mathcal{C}_\theta^{(n)}}P_nm(\theta,t).
\end{equation}
This  coincides with $~\widehat{c_\theta}$, the solution of the
system of equations (\ref{syst empirique}). So, we can write
\begin{equation}\label{c chapeau n duale}
    \widehat{c_\theta}=\arg\sup_{t\in\mathcal{C}_\theta^{(n)}} P_nm(\theta,t).
\end{equation}
Using (\ref{estim de phiMtheta dual}), we obtain the following
representations for the estimates
$\widehat{\phi}(\mathcal{M},P_0)$ in (\ref{estim de phiM  cas
cont}) and $\widehat{\theta}_\phi$ in (\ref{estim de theta0 cas
cont})
\begin{equation}\label{estim de phiM dual}
  \widehat{\phi}(\mathcal{M},P_0)=\inf_{\theta\in\Theta}
  \sup_{t\in\mathcal{C}_\theta^{(n)}}P_nm(\theta,t)
\end{equation}
and
\begin{equation}\label{estm theta0 dual}
    \widehat{\theta}_\phi=\arg\inf_{\theta\in\Theta}
    \sup_{t\in\mathcal{C}_\theta^{(n)}}P_nm(\theta,t),
\end{equation}
respectively.\\
\\
\\
Formula (\ref{sup P0 m}) also has the following basic interest :
Consider the function
\begin{equation}\label{P0 m theta t}
    t\in\mathcal{C}_\theta \mapsto P_0m(\theta,t),
\end{equation}
In order for integral (\ref{P0 m theta t}) to be properly defined,
we assume that
\begin{equation}\label{condition de quasi integrabilite}
\int \left|g_i(x,\theta)\right|~dP_0(x)<\infty,~\text{ for all
}~i\in\left\{1,\ldots,l\right\}.
\end{equation}
The domain of the function (\ref{P0 m theta t}) is
\begin{equation}\label{domaine de P0 m theta t}
    \mathcal{D}_\phi(\theta):=\left\{t\in\mathcal{C}_\theta ~
    \text{ such that }~
    P_0m(\theta,t)>-\infty\right\}.
\end{equation}
\noindent The function $t\mapsto P_0m(\theta,t)$ is strictly
concave
 on the convex set $\mathcal{D}_\phi(\theta)$. Whenever it
has a maximum $t^*$, then it is unique, and if it belongs to the
interior of $\mathcal{D}_\phi(\theta)$, then it satisfies the
first order condition. Therefore $t^*$ satisfies   system
(\ref{syst teta fixe}). In turn, this implies that the measure
$Q^*$ defined through
$dQ^*:={\varphi'}^{-1}\left({t^*}^T\overline{g}(\theta)\right)~dP_0$
is the projection of $P_0$ on $\Omega$, by Theorem 3.4  part 1 in
\cite{Bronia_Kez2006_STUDIA}. This implies that $Q^*$ and $P_0$
share the same support. We summarize the above arguments as
follows

\begin{proposition}\label{existence de Q star}
Assume that (\ref{condition de quasi integrabilite}) holds and
that
\begin{enumerate}
 \item[(i)] there exists some $s$ in the interior of $\mathcal{D}_\phi(\theta)$
  such that for all $t$ in $\partial \mathcal{D}_\phi(\theta)$, the frontier of
$\mathcal{D}_\phi(\theta)$, it holds $P_0m(\theta,t)<
P_0m(\theta,s)$;
 \item[(ii)] for all $t$ in
the interior of $\mathcal{D}_\phi(\theta)$, there exists
  a neighborhood $V(t)$ of $t$, such that the classes of functions
 $\left\{x\rightarrow \frac{\partial}{\partial r_i}m(x,\theta,r),~~r\in V(t) \right\}$
 are dominated  ($P_0$-a.s.) by some
 $P_0$-integrable function $x\rightarrow
H(x,\theta)$.
\end{enumerate}
Then $P_0$ admits an unique projection $Q_\theta^*$ on
$\mathcal{M}_\theta$ having the same support as $P_0$ and
\begin{equation}\label{form explicit de Q star par dualite}
    dQ_\theta^*=\varphi'^{-1}\left({c_\theta}^T\overline{g}(\theta)\right)dP_0,
\end{equation}
where $c_\theta$ is the unique solution of the system of equations
(\ref{syst teta fixe}).
\end{proposition}
\begin{remark}
In the case of $KL$-divergence, comparing this Proposition  with
Theorem 3.3 in \cite{Csiszar1975}, we observe that the dual
formula (\ref{sup P0 m}) provides weaker conditions on the class
of functions
$\left\{\overline{g}(\theta),~\theta\in\Theta\right\}$ than the
geometric approach.
\end{remark}

\begin{remark}\label{condition d'existence de Qstar cq}
The result of \cite{BorweinLewis1991}, with some  additional
conditions, provides more practical tools for obtaining the
results in Proposition \ref{existence de Q star}. Assume that the
functions $g_j(\theta)$ belongs to the space
$L_p(\mathcal{X},P_0)$ with $1\leq p\leq \infty$ and that the
following ``constraint qualification'' holds
 \begin{equation}\label{constraint qualification P}
\text{there exists some } Q_0 \text{ in } \mathcal{M}_\theta
\text{ such that : }~ a< \inf \frac{dQ_0}{dP_0}\leq \sup
\frac{dQ_0}{dP_0}<b,
 \end{equation}
 with $(a,b)$ is the domain $D_\varphi$ of the divergence function $\varphi$ and
 $\mathcal{M}_\theta$ is the set of all signed measures $Q$
 a.c. w.r.t. $P_0$, satisfying the linear constraints
 and  such that $\frac{dQ}{dP_0}$ belong to $L_q(\mathcal{X},P_0)$,
 ($1\leq q\leq \infty$ and $1/p+1/q=1$). In this case, applying
 Corollary 2 in \cite{BorweinLewis1991}, we obtain
 \begin{equation*}
   \phi\left(\mathcal{M}_\theta,P_0\right)=\sup_{t\in\mathbb{R}^{(l+1)}}
   \left\{t_0-\int \varphi^*\left(t^T\overline{g}(x,\theta)\right)~dP_0(x)
   \right\}
 \end{equation*}
 (with dual attainement). Furthermore, if for a dual optimal
 $c_\theta$, it holds
 \begin{equation*}
  \lim_{y\downarrow -\infty} \frac{\varphi(y)}{y} < \inf_x
  c_\theta^T\overline{g}(x,\theta)\leq
  \sup_x
  c_\theta^T\overline{g}(x,\theta) <
  \lim_{y\uparrow +\infty} \frac{\varphi(y)}{y}~\text{ for all
  }x ~(P_0~a.s.),
 \end{equation*}
 then the unique projection $Q^*_\theta$ of $P_0$ on
 $\mathcal{M}_\theta$
 is given by
\begin{equation}\label{form explicit de Q star par BL}
    dQ_\theta^*={\varphi^*}'\left({c_\theta}^T\overline{g}(\theta)\right)dP_0.
\end{equation}
Note that if $\varphi^*$ is strictly convex, then $c_\theta$ is
unique and
\begin{equation*}
\sup_{t\in\mathbb{R}^{(l+1)}}
   \left\{t_0-\int \varphi^*\left(t^T\overline{g}(x,\theta)\right)~dP_0(x)
   \right\}=\sup_{t\in\mathcal{C}_\theta}
   \left\{t_0-\int \varphi^*\left(t^T\overline{g}(x,\theta)\right)~dP_0(x)
   \right\}, \text{ and }
\end{equation*}
\begin{equation*}
 {\varphi^*}'\left({c_\theta}^T\overline{g}(\theta)\right)={\varphi'}^{-1}
 \left({c_\theta}^T\overline{g}(\theta)\right).
\end{equation*}
\cite{Leonard2001a} and \cite{Leonard2001b} gives, under minimal
conditions, duality theorems of minimum $\phi$-divergences  and
characterization of projections under
 linear constraints,
 which generalize the results
 given by \cite{BorweinLewis1991} and
 \cite{BorweinLewis1993}. These results are used recently
 by \cite{Bertail2003} and \cite{Bertail2006} in empirical likelihood.
\end{remark}

\section{Asymptotic properties and Statistical Tests}
\noindent In the sequel, we assume that the conditions in
Proposition \ref{existence de Q star n} (or  Remark
(\ref{condition d'existence de Qn cq})) and in Proposition
\ref{existence de Q star} (or Remark (\ref{condition d'existence
de Qstar cq})) hold. This allows to use the representations
(\ref{estim de phiMtheta dual}), (\ref{estim de phiM dual}) and
(\ref{estm theta0 dual}) in order to study the asymptotic behavior
of the proposed estimates (\ref{estim de phiMtheta cas cont}),
(\ref{estim de phiM cas cont}) and (\ref{estim de theta0 cas
cont}). All the results in the present Section are obtained
through classical methods of parametric statistics; see e.g.
\cite{vanderVaart1998} and \cite{Sen-Singer1993}.
We first consider the case when $\theta$ is fixed, and we study
the asymptotic behavior of  the estimate
$\widehat{\phi}(\mathcal{M}_{\theta},P_0)$ (see (\ref{estim de
phiMtheta cas cont})) of
$\phi(\mathcal{M}_{\theta},P_0):=\inf_{Q\in \mathcal{M}_{\theta
}}\phi(Q,P_0)$ both when $P_0\in\mathcal{M}_\theta$ and when
$P_0\not\in\mathcal{M}_\theta$. This is done in the first
Subsection. In the second Subsection, we study the asymptotic
behavior of the EM$\phi$D estimates $\widehat{\theta}_\phi$ and
the estimates $\widehat{\phi}\left(\mathcal{M},P_0\right)$ both in
the two cases when $P_0$ belongs to $\mathcal{M}$ and when $P_0$
does not belong to $\mathcal{M}$. The solution of \textit{Problem
1} is given in Subsection 3.3 while \textit{Problem 2} is treated
in Subsections 3.1, 3.2, 3.3 and 3.4.

\subsection{Asymptotic properties of the estimates for a given
$\theta\in\Theta$ } First we state consistency.
\subsection*{Consistency}
We state both weak and strong consistency of the estimates
$\widehat{c_\theta}$ and $\widehat{\phi}(\mathcal{M}_\theta,P_0)$
using their representations (\ref{c chapeau n duale}) and
(\ref{estim de phiMtheta dual}), respectively. Denote $\|.\|$ the
Euclidian norm defined on $\mathbb{R}^d$ or on $\mathbb{R}^{l+1}$.
In order to state consistency, we need to define
\begin{equation*}
T_\theta := \left\{t\in\mathcal{C}_\theta ~\text{ such that }~
P_0m(\theta,t) > - \infty\right\},
\end{equation*}
and denote $T_\theta^c$ the complementary of the set $T_\theta$ in
the set $\mathcal{C}_\theta $, namely
\begin{equation*}
T_\theta^c := \left\{t\in\mathcal{C}_\theta ~\text{ such that }~
P_0m(\theta,t)= -\infty\right\}.
\end{equation*}
Note that, by Proposition \ref{existence de Q star},  the set
$T_\theta$ contains $c_\theta$.\\

\noindent We will consider the following condition
\begin{enumerate}
 \item [(C.1)] $\sup_{t\in
T_\theta}\left|P_nm(\theta,t)-P_0m(\theta,t)\right|$ converges to
$0$ a.s. (resp. in probability);
 \item [(C.2)] there exists $M<0$
and $n_0>0$, such that, for all $n>n_0$, it holds $\sup_{t\in
T_\theta^c}P_nm(\theta,t)\leq M$~
a.s. (resp. in probability).\\
\end{enumerate}

\noindent The condition  (C.2) makes sense, since for all $t\in
T_\theta^c$ we have $P_0m(\theta,t)=-\infty$.\\

\noindent Since the function $t\in T_\theta\mapsto P_0m(\theta,t)$
is strictly concave, the maximum $c_\theta$ is isolated, that is
\begin{equation}\label{c  theta est bien separer}
\text{ for any positive }~~  \epsilon,~ \text{ we have }~
\sup_{\{t\in\mathcal{C}_\theta~:~\|t-c_\theta\|\geq\epsilon\}}
P_0m(\theta,t)<P_0m(\theta,c_\theta).
\end{equation}

\begin{proposition}\label{consitence de l estimateur de ctheta}
Assume that conditions (C.1) and (C.2)  hold. Then
\begin{enumerate}
 \item [(i)]  the estimates
$\widehat{\phi}(\mathcal{M}_\theta,P_0)$ converge to
$\phi(\mathcal{M}_\theta,P_0)$ a.s. (resp. in probability).
 \item [(ii)]
the estimates $\widehat{c_\theta}$ converge to $c_\theta$ a.s.
(resp. in probability).
\end{enumerate}
\end{proposition}

\subsection*{Asymptotic distributions}
Denote $~m'(\theta,t)~$ the $~(l+1)$-dimensional vector with
entries $~\frac{\partial}{\partial t_i}m(\theta,t)$,
$~m''(\theta,t)~$ the $~(l+1)\times (l+1)$-matrix with entries
$~\frac{\partial^2}{\partial t_i \partial t_j}m(\theta,t)$,
$~\underline{0}_l:=(0,\ldots,0)^T\in\mathbb{R}^l$,
$~\underline{0}_d:=(0,\ldots,0)^T\in\mathbb{R}^d$,
$~\underline{c}$ the $~(l+1)-$vector defined by
$~\underline{c}:={\left(0,\underline{0}_l^T\right)}^T$, and
$~P_0g(\theta)g(\theta)^T~$ the $~l\times l-$matrix defined by
\begin{equation*}
P_0g(\theta)g(\theta)^T:=
{\left[P_0g_i(\theta)g_j(\theta)\right]}_{i,j=1,\ldots,l}.
\end{equation*}

We will consider the following assumptions
\begin{enumerate}
 \item [(A.1)] $~\widehat{c_\theta}~$ converges in probability to $~c_\theta$;
 \item [(A.2)] the function $~t\mapsto m(x,\theta,t)~$ is
 $~\mathcal{C}^3~$ on a neighborhood $~V(c_\theta)~$ of $~c_\theta~$ for all $x$
  ~($P_0$-a.s.), and
 all partial derivatives of order $~3~$ of the function
 $\left\{t\mapsto m(x,\theta,t),~t\in V(c_\theta)\right\}$ are dominated  by some
 $~P_0$-integrable function $~x\mapsto H(x)$;
 \item [(A.3)] $~P_0\left(\|m'(\theta,c_\theta)\|^2\right)~$ is finite, and
 the matrix $~P_0m''(\theta,c_\theta)~$ exists and is invertible.
\end{enumerate}

\begin{theorem}\label{lois limite a theta fixe}
Assume that assumptions (A.1-3) hold. Then
 \begin{enumerate}
  \item [(1)] $\sqrt{n}\left(\widehat{c_\theta}-c_\theta\right)$ converges to a
  centered normal multivariate variable with covariance matrix
  \begin{equation}
   V={\left[-P_0m''(\theta,c_\theta)\right]}^{-1}
   \left[P_0m'(\theta,c_\theta)m'(\theta,c_\theta)^T\right]
   {\left[-P_0m''(\theta,c_\theta)\right]}^{-1}.
  \end{equation}
In the special case, when  $~P_0~$ belongs to
$\mathcal{M}_\theta$, then
  $c_\theta=\underline{c}$ and
  \begin{equation}\label{variance si P est dans Mtheta}
   V=\varphi''(1)^2\left[
\begin{array}{cc}
  0 &  \underline{0}_l^T\\
  \underline{0}_l & {\left[P_0g(\theta)g(\theta)^T\right]}^{-1} \\
\end{array}
\right].
  \end{equation}
 \item [(2)] If $P_0$ belongs to $\mathcal{M}_\theta$, then  the
 statistics
\begin{equation*}
\frac{2n}{\varphi''(1)}\widehat{\phi}
 \left(\mathcal{M}_\theta,P_0\right)
\end{equation*}
converge in distribution to a $\chi^2$ variable with $l$ degrees
of freedom.
 \item [(3)] If $P_0$ does not belong to $\mathcal{M}_\theta$, then
\begin{equation*}
\sqrt{n}\left(\widehat{\phi}
 \left(\mathcal{M}_\theta,P_0\right)-\phi(\mathcal{M}_\theta,P_0)\right)
\end{equation*}
converges to a centered normal variable with variance
\begin{equation*}
\sigma^2:=P_0m(\theta,c_\theta)^2-\left(P_0m(\theta,c_\theta)\right)^2.\\
\end{equation*}
\end{enumerate}
\end{theorem}

\vskip 1cm

\begin{remark}
\begin{enumerate}
 \item [(a)] When specialized to the modified Kullback-Leibler
divergence, Theorem \ref{lois limite a theta fixe} part (2) gives
the limiting distribution of the empirical log-likelihood ratio
$2n\widehat{KL_m}(\mathcal{M}_\theta,P_0)$ which is the result in
\cite{Owen1990} Theorem 1. Part (3) gives its limiting
distribution when $P_0$ does not belong to $\mathcal{M}_\theta$.
 \item[(b)] Nonparametric confidence regions  $(CR_\phi)$ for  $\theta_0$
 of asymptotic level $(1-\epsilon)$ can be constructed using the
statistics
\begin{equation*}
\frac{2n}{\varphi''(1)}\widehat{\phi}
 \left(\mathcal{M}_\theta,P_0\right),
\end{equation*}
through
\begin{equation*}
    CR_\phi:=\left\{\theta\in\Theta~\text{ such that }~\frac{2n}{\varphi''(1)}\widehat{\phi}
 \left(\mathcal{M}_\theta,P_0\right)\leq
 ~q_{(1-\epsilon)}\right\},
\end{equation*}
where $(1-\epsilon)$ is the  $(1-\epsilon)$-quantile of a
$\chi^2(l)$ distribution. It would be interesting to obtain the
divergence leading to optimal confidence regions in the sense of
\cite{Neyman1937} (see \cite{Takagi1998}), or the optimal
divergence  leading to confidence regions with  small length
(volume, area or diameter) and covering the true value $\theta_0$
with large enough probability.
\end{enumerate}
\end{remark}

\subsection{Asymptotic properties of the estimates
$\widehat{\theta}_\phi$ and $\widehat{\phi}(\mathcal{M},P_0)$}
First we state consistency.
\subsection*{Consistency}
We assume that when $P_0$ does not  belong to the model
$\mathcal{M}$, the minimum, say $\theta^*$, of the function
$\theta\in\Theta \mapsto \inf_{Q\in\mathcal{M}_\theta}\phi(Q,P_0)$
exists and  is unique. Hence $P_0$ admits a projection  on
$\mathcal{M}$ which we denote $Q^*_{\theta^*}$. Obviously when
$P_0$ belongs to the model $\mathcal{M}$, then $\theta^*=\theta_0$
and $Q^*_{\theta^*}=P_0$. We will consider the following
conditions
\begin{enumerate}
\item [(C.3)] $\sup_{\{\theta\in\Theta, t\in T_\theta\}} \left|
P_nm(\theta,t)-P_0m(\theta,t)\right|$ tends to $0$ a.s. (resp. in
probability); \item [(C.4)] there exists a neighborhood
$V(c_{\theta^*})$ of $c_{\theta^*}$ such that
\begin{enumerate}
\item[(a)] for any positive $\epsilon$, there exists some positive
$\eta$ such that for all $t\in V(c_{\theta^*})$ and all
$\theta\in\Theta$ satisfying $\|\theta-\theta^*\|\geq\epsilon$, it
holds $P_0m(\theta^*,t)<P_0m(\theta,t)-\eta$;

\item [(b)] there exists some function $H$ such that for all $t$
in $V(c_{\theta^*})$, ~we have $\left|m(t,\theta_0)\right|\leq
H(x)$~ ($P_0$-a.s.)~with $P_0H < \infty$;
\end{enumerate}
\item[(C.5)] there exits $M<0$ and $n_0>0$ such that for all
$n\geq n_0$, we have
\begin{equation}
\sup_{\theta\in\Theta} \sup_{t\in T_\theta^c}P_n m(\theta,t)\leq M
~ \text{ a.s. (resp. in probability).}
\end{equation}
\end{enumerate}

\begin{proposition}\label{consistence de theta n et}
Assume that conditions (C.3-5) hold. Then
\begin{enumerate}
\item [(i)] the estimates $\widehat{\phi}(\mathcal{M},P_0)$
converge  to $\phi(\mathcal{M},P_0)$  a.s. (resp. in probability).
\item[(ii)]
$\sup_{\theta\in\Theta}\left\|\widehat{c_\theta}-c_{\theta}\right\|$
converge to $0$ a.s. (resp. in probability).

\item[(iii)] The ME$\phi$D estimates $\widehat{\theta}_\phi$
converge to $\theta^*$ a.s. (resp. in probability).
\end{enumerate}
\end{proposition}

\subsection*{Asymptotic distributions}
When $P_0\in\mathcal{M}$, then by assumption, there exists unique
$\theta_0\in\Theta$ such that $P_0\in\mathcal{M}_{\theta_0}$.
Hence $\theta^*=\theta_0$ and
$c_{\theta^*}=c_{\theta_0}=\underline{c}$. We state the limit
distributions of the estimates $\widehat{\theta}_\phi$ and
    $\widehat{c_{\widehat{\theta}_\phi}}$ when $P_0\in\mathcal{M}$
    and when $P_0\not\in\mathcal{M}$. We will make use of the following
assumptions
\begin{enumerate}
  \item [(A.4)] Both estimates $\widehat{\theta}_\phi$ and
    $\widehat{c_{\widehat{\theta}_\phi}}$ converge in probability respectively to
    $\theta^*$ and $c_{\theta^*}$;
 \item [(A.5)] the function $(\theta,t)\mapsto m(x,\theta,t)$ is
 $\mathcal{C}^3$ on some neighborhood $V(\theta^*,c_{\theta^*})$
 for all $x$ ~($P_0$-a.s.), and the partial derivatives of order $3$ of the
 functions\\ $\left\{(\theta,t)\mapsto m(x,\theta,t),~
 (\theta,t)\in V(\theta^*,c_{\theta^*})\right\}$ are dominated
  by some $P_0-$integrable function $H(x)$;
 \item [(A.6)] $P_0\left(\left\|\frac{\partial}{\partial t}m(\theta^*,
 c_{\theta^*})\right\|^2\right)$ and $P_0\left(\left\|\frac{\partial}{\partial
 \theta}m(\theta^*,c_{\theta^*})\right\|^2\right)$ are finite, and the
 matrix
\begin{equation*}
S:=\left(%
\begin{array}{cc}
  S_{11} & S_{12} \\
  S_{21} & S_{22} \\
\end{array}%
\right),
\end{equation*}
with $S_{11}:=P_0\frac{\partial^2}{\partial t^2}m(\theta^*,
 c_{\theta^*})$, $S_{12}={S_{21}}^T:=P_0\frac{\partial^2}
 {\partial t\partial\theta}m(\theta^*,c_{\theta^*})$ and
$S_{22}:=P_0\frac{\partial^2}{\partial \theta^2}m(\theta^*,
 c_{\theta^*})$,
 exists and is invertible.
\end{enumerate}

\begin{theorem}\label{loi limite de thetan et cn}
Let $P_0$ belongs to $\mathcal{M}$ and assumptions (A.4-6) hold.
Then, both $\sqrt{n}\left(\widehat{\theta}_\phi-\theta_0\right)$
and
$\sqrt{n}\left(\widehat{c_{\widehat{\theta}_\phi}}-\underline{c}\right)$
converge in distribution to a centered  multivariate normal
variable with covariance matrix, respectively
\begin{equation}\label{LCM de thetan}
    V={\left\{\left[P_0\frac{\partial}{\partial\theta}g(\theta_0)\right]
  \left[P_0\left(g(\theta_0)g(\theta_0)^T\right)\right]^{-1}\left[P_0\frac{\partial
  }{\partial\theta}g(\theta_0)\right]^T\right\}}^{-1},
\end{equation}
and
\begin{eqnarray}\label{CM 11}
    U & = & \varphi''(1)^2\left[
\begin{array}{cc}
  0 & \underline{0}_l^T \\
  \underline{0}_l & \left[P_0g(\theta_0)g(\theta_0)^T\right]^{-1} \\
\end{array}
\right]-\varphi''(1)^2\left[
\begin{array}{cc}
  0 & \underline{0}_l^T \\
  \underline{0}_l & \left[P_0g(\theta_0)g(\theta_0)^T\right]^{-1} \\
\end{array}
\right]\times\nonumber\\
& &
\times\left[\underline{0}_d,P_0\frac{\partial}{\partial\theta}g(\theta_0)\right]^T
V
\left[\underline{0}_d,P_0\frac{\partial}{\partial\theta}g(\theta_0)\right]\left[
\begin{array}{cc}
  0 & \underline{0}_l^T \\
  \underline{0}_l & \left[P_0g(\theta_0)g(\theta_0)^T\right]^{-1} \\
\end{array}
\right], \nonumber
\end{eqnarray}
and the estimates $\widehat{\theta}_\phi$ and
$\widehat{c_{\widehat{\theta}_\phi}}$ are asymptotically
uncorrelated.\\
\end{theorem}

\begin{remark}\label{tous les estim sont equiv}
When specialized to the modified Kullback-Leibler divergence, the
estimate $\widehat{\theta}_{KL_m}$ is the empirical likelihood
estimate (ELE) (noted $\widetilde{\theta}$ in
\cite{Qin-Lawless1994}), and the above result gives the limiting
distribution of $\sqrt{n}(\widehat{\theta}_{KL_m}-\theta_0)$ which
coincides with the result in Theorem 1 in \cite{Qin-Lawless1994}.
Note also that all ME$\phi$DE's including ELE have the same
limiting distribution with the same variance when $P_0$ belongs to
$\mathcal{M}$. Hence they are all equally first order efficient.
\end{remark}

\begin{theorem}\label{loi limite du vect thetan  cn}
Assume that $P_0$ does not belong to $\mathcal{M}$ and that
assumptions (A.4-6) hold. Then
\begin{equation*}
    \sqrt{n}\left(%
\begin{array}{c}
  \widehat{c_{\widehat{\theta}_\phi}}-c_{\theta^*} \\
  \widehat{\theta}_\phi-\theta^* \\
\end{array}%
\right)
\end{equation*}
converges in distribution to a centered multivariate normal
variable with covariance matrix
\begin{equation*}\label{matrice limit a l'exter du model}
  W = S^{-1}M{S^{-1}}
\end{equation*}
where
\begin{equation*}
M := P_0\left(\left(%
\begin{array}{c}
  \frac{\partial}{\partial t}m\left(\theta^*,c_{\theta^*}\right) \\
  \frac{\partial }{\partial \theta}m\left(\theta^*,c_{\theta^*}\right) \\
\end{array}%
\right)\left(%
\begin{array}{c}
  \frac{\partial}{\partial t}m\left(\theta^*,c_{\theta^*}\right) \\
  \frac{\partial }{\partial \theta}m\left(\theta^*,c_{\theta^*}\right) \\
\end{array}%
\right)^T\right).
\end{equation*}
$\theta^*$ and $c_{\theta^*}$ are characterized by
\begin{equation*}
\theta^*:=\arg\inf_{\theta\in\Theta}\phi\left(\mathcal{M}_\theta,P_0\right),
\end{equation*}
\begin{equation*}
dQ^*_{\theta^*}={\varphi'}^{-1}\left(c_{\theta^*}^T\overline{g}(\theta)\right)
dP_0 ~\text{ and }~ Q^*_{\theta^*}\in\mathcal{M}_{\theta^*}.
\end{equation*}
\end{theorem}

\subsection{Tests of model}

In order to test the hypothesis $\mathcal{H}_0$~:~$P_0$ belongs to
$\mathcal{M}$ against the alternative $\mathcal{H}_1$~:~$P_0$ does
not belong to $\mathcal{M}$, we can use the estimates
$\widehat{\phi}(\mathcal{M},P_0)$ of $\phi(\mathcal{M},P_0)$, the
$\phi-$divergences between the model $\mathcal{M}$ and the
distribution $P_0$. Since $\phi(\mathcal{M},P_0)$ is nonnegative
and take value $0$ only when $P_0$ belongs to $\mathcal{M}$
(provided that $P_0$ admits a projection on $\mathcal{M}$), we
reject the hypothesis $\mathcal{H}_0$ when the estimates take
large values. In the following Corollary, we give the asymptotic
law of the estimates $\widehat{\phi}(\mathcal{M},P_0)$ both under
$\mathcal{H}_0$ and under $\mathcal{H}_1$.

\begin{corollary}\label{the pour le test du model}\
\begin{enumerate}
 \item[(i)] Assume that the assumptions of Theorem \ref{loi limite de
thetan et cn} hold and that
  $l>d$. Then, under $\mathcal{H}_0$, the statistics
 \begin{equation*}
    \frac{2n}{\varphi''(1)}\widehat{\phi}(\mathcal{M},P_0)
\end{equation*}
converge in distribution to a $\chi^2$ variable with $(l-d)$
degrees of freedom.
 \item[(ii)] Assume that the assumptions of
Theorem \ref{loi limite du vect thetan  cn} hold. Then, under
$\mathcal{H}_1$,~ we have :
\begin{equation}
\sqrt{n}\left(\widehat{\phi}
  (\mathcal{M},P_0)-\phi(\mathcal{M},P_0)\right)
\end{equation}
converges to centered normal variable with variance
\begin{equation*}
\sigma^2=P_0m(\theta^*,c_{\theta^*})^2-
  \left(P_0m(\theta^*,c_{\theta^*})\right)^2
\end{equation*}
where $\theta^*$ and $c_{\theta^*}$ satisfy
\begin{equation*}
\theta^*:=\arg\inf_{\theta\in\Theta}\phi\left(\mathcal{M}_\theta,P_0\right),
\end{equation*}
\begin{equation*}
    \varphi'\left(\frac{dQ_{\theta^*}^*}{dP_0}(x)\right)=
    c_{\theta^*}^T\overline{g}(x,\theta^*)
    ~~~~\text{ and }~~~~ Q_{\theta^*}^*\in\mathcal{M}_{\theta^*}.
\end{equation*}
\end{enumerate}
\end{corollary}

\begin{remark}
This Theorem allows to perform tests of model  of asymptotic level
$\alpha$; the critical regions are
\begin{equation}\label{region critique pr le modele}
    C_\phi:=\left\{\frac{2n}{\varphi''(1)}
    \widehat{\phi}(\mathcal{M},P_0)>q_{(1-\alpha)}\right\},
\end{equation}
where $q_{(1-\alpha)}$ is the $(1-\alpha)-$quantile of the
$\chi^2$ distribution with $(l-d)$ degrees of freedom. Also these
tests are all asymptotically powerful, since the estimates
$\widehat{\phi}(\mathcal{M},P_0)$ are $n-$consistent estimates of
$\phi(\mathcal{M},P_0)=0$ under $\mathcal{H}_0$ and
$\sqrt{n}-$consistent estimates of $\phi(\mathcal{M},P_0)$ under
$\mathcal{H}_1$.\\
\end{remark}

\noindent We assume now that the p.m. $P_0$ belongs to
$\mathcal{M}$. We will perform simple and composite tests on the
parameter $\theta_0$ taking  into account of the information
$P_0\in\mathcal{M}$.

\subsection{Simple tests on the parameter}
Let
\begin{equation}\label{simple test sur param}
    \mathcal{H}_0 ~~:~~\theta_0=\theta_1~~~\text{ versus }~~~
    \mathcal{H}_1 ~~:~~\theta_0\in\Theta\setminus\{\theta_1\},
\end{equation}
where $\theta_1$  is a given known value. We can use the following
statistics to perform tests pertaining to (\ref{simple test sur
param})
\begin{equation*}
S_n^\phi:=\widehat{\phi}(\mathcal{M}_{\theta_1},P_0)-\inf_{\theta\in\Theta}
\widehat{\phi}(\mathcal{M}_\theta,P_0).
\end{equation*}
Since
\begin{equation*}
\phi(\mathcal{M}_{\theta_1},P_0)-\inf_{\theta\in\Theta}
\phi(\mathcal{M}_\theta,P_0)=\phi(\mathcal{M}_{\theta_1},P_0)
\end{equation*}
are nonnegative and take value $0$ only when $\theta_0=\theta_1$,
we reject the hypothesis $\mathcal{H}_0$ when the statistics
$S_n^\phi$ take large values.\\

\noindent We give the limit distributions of the statistics
$S_n^\phi$ in the following Corollary which we can prove using
some algebra and
 arguments used in the proof of Theorem \ref{loi limite de thetan
et cn} and Theorem \ref{loi limite du vect thetan  cn}.

\begin{corollary}\label{tests simple sur le param}\
\begin{enumerate}
 \item[(i)] Assume that assumptions of Theorem \ref{loi limite de thetan et
cn} hold. Then under $\mathcal{H}_0$, the statistics
\begin{equation*}
\frac{2n}{\varphi''(1)}S_n^\phi
\end{equation*}
converge in distribution to $\chi^2$ variable with $d$ degrees of
freedom.
 \item[(ii)] Assume that assumptions of Theorem \ref{loi limite de thetan et
cn} hold. Then under $\mathcal{H}_1$,
\begin{equation*}
 \sqrt{n}\left(S_n^\phi-\phi\left(\mathcal{M}_{\theta_1},P_0\right)\right)
\end{equation*}
converges to a centered normal variable with variance
\begin{equation*}
 \sigma^2=P_0m(\theta_1,c_{\theta_1})^2-\left(P_0m(\theta_1,c_{\theta_1})\right)^2.
\end{equation*}
\end{enumerate}
\end{corollary}

\begin{remark}
When specialized to the $KL_m$-divergence, the statistic
$2nS_n^{KL_m}$ is the empirical likelihood ratio statistic (see
           \cite{Qin-Lawless1994} Theorem 2).
\end{remark}
\subsection{Composite tests on the parameter}

Let
\begin{equation}\label{h nonparam}
\begin{array}{ccccc}
h & : & \mathbb{R}^d & \rightarrow & \mathbb{R}^k
\end{array}
\end{equation}
be some function such that the $(d\times k)-$matrix
$H(\theta):=\frac{\partial}{\partial\theta}h(\theta)$ exists, is
continuous and has rank $k$ with $0<k<d$. Let us define the
composite null hypothesis
\begin{equation}\label{Theta 0}
    \Theta_0:=\left\{\theta\in\Theta ~\text{ such that }~
    h(\theta)=0\right\}.
\end{equation}
We consider the composite test
\begin{equation}\label{composite test nonparam}
    \mathcal{H}_0 ~~:~~\theta_0\in\Theta_0~~~\text{ versus }~~~
    \mathcal{H}_1 ~~:~~\theta_0\in\Theta\setminus\Theta_0,
\end{equation}
i.e., the test
\begin{equation}\label{composite test nonparam}
    \mathcal{H}_0 ~~:~~P_0\in \bigcup_{\theta\in\Theta_0}\mathcal{M}_\theta
    ~~~\text{ versus }~~~
    \mathcal{H}_1 ~~:~~P_0\in
    \bigcup_{\theta\in\Theta\setminus\Theta_0}\mathcal{M}_\theta.
\end{equation}
This test is equivalent to the following one
\begin{equation}\label{composite test nonparam}
    \mathcal{H}_0 ~~:~~\theta_0\in  f(B_0)~~~\text{ versus }~~~
    \mathcal{H}_1 ~~:~~\theta_0\not\in f(B_0),
\end{equation}
where $f~:~\mathbb{R}^{(d-k)}~\rightarrow ~\mathbb{R}^d$ is a
function such that the matrix
$G(\beta):=\frac{\partial}{\partial\beta}g(\beta)$ exists and has
rank $(d-k)$, and $B_0:=\left\{\beta\in\mathbb{R}^{(d-k)} \text{
such that } f(\beta)\in \Theta_0\right\}$. Therefore
$\theta_0\in\Theta_0$ is an equivalent  statement for
$\theta_0=f(\beta_0), \beta_0\in B_0$.\\

\noindent The following statistics are used to perform tests
pertaining to (\ref{composite test nonparam}) :
\begin{equation*}
T_n^\phi:=\inf_{\beta\in
B_0}\widehat{\phi}\left(\mathcal{M}_{f(\beta)},P_0\right)-
\inf_{\theta\in\Theta}\widehat{\phi}\left(\mathcal{M}_\theta,P_0\right).
\end{equation*}
Since
\begin{equation*}
\inf_{\beta\in
B_0}\phi(\mathcal{M}_{f(\beta)},P_0)-\inf_{\theta\in\Theta}
\phi(\mathcal{M}_\theta,P_0)=\inf_{\beta\in
B_0}\phi(\mathcal{M}_{f(\beta)},P_0)
\end{equation*}
are nonnegative and take value $0$ only when $\mathcal{H}_0$
holds, we reject the hypothesis $\mathcal{H}_0$ when the
statistics $T_n^\phi$ take large values.\\

\noindent We give the limit distributions of the statistics
$T_n^\phi$ in the following Corollary.

\begin{corollary}\label{the pour le test compose sur le param}\
 \begin{enumerate}
  \item[(i)] Assume that  assumptions of Theorem \ref{loi limite de
thetan et cn} hold. Under $\mathcal{H}_0$, the statistics
$T_n^\phi$ converge in distribution to a $\chi^2$ variable with
$(d-k)$ degrees of freedom.
  \item[(ii)] Assume that there exists $\beta^*\in B_0$, such
  that $\beta^*=\arg\inf_{\beta\in
  B_0}\phi\left(\mathcal{M}_{f(\beta)},P_0\right)$. If the
  assumptions of Theorem \ref{loi limite du vect thetan  cn} hold
  for $\theta^*=f(\beta^*)$, then
 \begin{equation*}
\sqrt{n}\left(T_n^\phi-\phi\left(\mathcal{M}_{\theta^*},P_0\right)\right)
 \end{equation*}
converges to a centered normal variable with variance
\begin{equation*}
 \sigma^2=P_0m(\theta^*,c_{\theta^*})^2-\left(P_0m(\theta^*,c_{\theta^*})\right)^2.
\end{equation*}
\end{enumerate}
\end{corollary}

\section{ Estimates of the distribution function through projected distributions}

\noindent In this Subsection, the measurable space
$\left(\mathcal{X},\mathcal{B}\right)$ is
$\left(\mathbb{R},\mathcal{B}_\mathbb{R}\right)$. For all
$\phi-$divergence, by (\ref{estim de phiM  cas cont}), we have
\begin{equation*}
    \widehat{\phi}\left(\mathcal{M},P_0\right)=\phi\left(\mathcal{M},P_n\right)=
    \phi\left(\widehat{Q^*_{\widehat{\theta}_\phi}},P_n\right).
\end{equation*}
Proposition \ref{Q satr n est atteint} above provides the
description of $\widehat{Q^*_{\widehat{\theta}_\phi}}$.\\

\noindent So, for all $\phi$-divergence, we  estimate  the
distribution function $F$ using
$\widehat{Q^*_{\widehat{\theta}_\phi}}$ the $\phi-$projection of
$P_n$ on $\mathcal{M}$, through
\begin{eqnarray}\label{F chapeau n}
    \widehat{F}_n(x) & := &
    \sum_{i=1}^{n}\widehat{Q^*_{\widehat{\theta}_\phi}}(X_i)
    \mathds{1}_{(-\infty,x]}(X_i)\nonumber\\
     & = &
     \frac{1}{n}\sum_{i=1}^n\overleftarrow{\varphi'}
     \left(\widehat{c_{\widehat{\theta}_\phi}}^T\overline{g}(X_i,
     \widehat{\theta}_\phi)\right)\mathds{1}_{(-\infty,x]}(X_i).
\end{eqnarray}

\begin{remark} When the estimating equation
\begin{equation}\label{estimating equation}
   \frac{1}{n}\sum_{i=1}^{n}g(X_i,\theta)=\underline{0}_d
\end{equation}
admits a solution $\widetilde{\theta}_n$, then $P_n$ belongs to
$\mathcal{M}$. If the solution is unique then
$\widehat{\theta}_\phi=\widetilde{\theta}_n$. Hence by Proposition
\ref{existence de Q star n}
\begin{equation*}
\text{ for all }~ i\in\left\{1,2,\ldots,n\right\},~ \text{ we have
}~~ \widehat{Q^*_{\widehat{\theta}_\phi}}(X_i)=\frac{1}{n},
\end{equation*}
and  $\widehat{F}_n(x)$, in this case, is the empirical cumulative
distribution function, i.e.,
\begin{equation*}
   \widehat{F}_n(x)=F_n(x):=\frac{1}{n}\sum_{i=1}^n\mathds{1}_{(-\infty,x]}(X_i).
\end{equation*}
So, the main interest is in the case where (\ref{estimating
equation}) does not admit a solution, that is in general when
$l>d$.
\end{remark}
\begin{remark}
The $\phi$-projections $\widehat{Q^*_{\widehat{\theta}_\phi}}$ of
$P_n$ on $\mathcal{M}$ may be  signed measures. For all
$\phi$-divergence satisfying $D_\varphi=\mathbb{R}_+^*$, the
$\phi$-projection $\widehat{Q^*_{\widehat{\theta}_\phi}}$ is a
p.m. if it exists. (for example, $KL_m$, $KL$, Hellinger, and
$\chi^2_m$ divergences all provide  p.m.'s).\\
\end{remark}

\noindent We give the limit law of the estimates $\widehat{F}_n$
of the distribution function $F$ in the following Theorem. We will
see that the estimate $\widehat{F}_n(x)$ is generally more
efficient than the empirical cumulative distribution function
$F_n(x)$.

\begin{theorem}\label{loi limite de F chapeau n}
Under the assumptions of Theorem \ref{loi limite de thetan et cn},
 $\sqrt{n}\left(\widehat{F}_n(x)-F(x)\right)$ converges in
distribution to a centered normal variable with variance
\begin{equation}\label{variance de F n chapeau}
    W(x)=F(x)\left(1-F(x)\right)-
    {\left[P_0\left(g(\theta_0)\mathds{1}_{(-\infty,x]}\right)\right]}^T
    \Gamma
    \left[P_0\left(g(\theta_0)\mathds{1}_{(-\infty,x]}\right)\right],
\end{equation}
with
\begin{eqnarray}
 \Gamma & = & {\left[P_0g(\theta_0)g(\theta_0)^T\right]}^{-1}
 - {\left[P_0g(\theta_0)g(\theta_0)^T\right]}^{-1}{\left[P_0
 \frac{\partial}{\partial\theta}g(\theta_0)\right]}^T V\times\nonumber\\
     & & \times \left[P_0\frac{\partial}{\partial\theta}g(\theta_0)\right]
 {\left[P_0g(\theta_0)g(\theta_0)^T\right]}^{-1},\nonumber
\end{eqnarray}
and
\begin{equation*}
   V={\left\{\left[P_0\frac{\partial}{\partial\theta}g(\theta_0)\right]
  \left[P_0\left(g(\theta_0)g(\theta_0)^T\right)\right]^{-1}\left[P_0\frac{\partial
  }{\partial\theta}g(\theta_0)\right]^T\right\}}^{-1}.
\end{equation*}
\end{theorem}

\section{Empirical likelihood and related methods}

\noindent In the present setting, the empirical likelihood (EL)
approach for the estimation of the parameter $\theta_{0}$ can be
summarized as follows. For any $\theta$ in $\Theta$, define the
\textit{profile likelihood ratio} of the sample
$X:=(X_{1},...,X_{n})$ through
\begin{equation*}
L_{n}(\theta ):=\sup \left\{ \prod\limits_{i=1}^{n}nQ(X_i)\text{
where }  Q(X_i)\geq 0,~\sum_{i=1}^n
Q(X_i)=1,~\sum_{i=1}^ng(X_i,\theta)Q(X_i)=0\right\} .
\end{equation*}
The estimate of $\theta _{0}$ through empirical likelihood (EL)
approach is then defined by
\begin{equation}
\widehat{\theta}_{EL}:=\arg \sup_{\theta \in \Theta
}L_{n}(\theta). \label{EL}
\end{equation}
The paper by \cite{Qin-Lawless1994} introduces
$\widehat{\theta}_{EL}$ and presents its properties. In this
Section, we show that $\widehat{\theta}_{EL}$ belongs to the
family of ME$\phi$D estimates for the specific choice
$\varphi(x)=-\log x+x-1$. We also discuss
 the problem of the existence of the solution of $\left(
\ref{EL}\right)$ for all $n$.\\

\noindent When $\varphi(x)=-\log x+x-1$, formula (\ref{estim de
theta0 cas cont}) clearly coincides with $\widehat{\theta}_{EL}$.
For test of hypotheses given by
$\mathcal{H}_0~:~P_0\in\mathcal{M}_\theta$ ~against~
$\mathcal{H}_1~:~P_0\not\in\mathcal{M}_\theta$ or for construction
of nonparametric confidence regions for $\theta_0$, the statistic
$2n \widehat{KL_m}(\mathcal{M}_\theta,P_0)$ coincides with the
empirical log-likelihood ratio  introduced  in \cite{Owen1988},
\cite{Owen1990} and  \cite{Qin-Lawless1994}. We state the results
of Section 3 in the present context. We will see that the approach
of empirical likelihood by divergence minimization, using the dual
representation of the $KL_m$-divergence and the explicit form of
the $KL_m$-projection of $P_0$, yields to the limit distribution
of the  statistic $2n \widehat{KL_m}(\mathcal{M}_\theta,P_0)$
under $\mathcal{H}_1$, which can not be achieved using the
approach in \cite{Owen1990} and \cite{Qin-Lawless1994}. Consider
\begin{equation*}
\widehat{\theta}_{KL_m}=\arg\inf_{\theta\in\Theta}\widehat{KL_m}(\mathcal{M}_\theta,P_0)
\end{equation*}
where
\begin{equation}\label{klm M theta}
    \widehat{KL_m}(\mathcal{M}_\theta,P_0)=\sup_{t\in\mathcal{C}_\theta}
    P_nm(\theta,t)
\end{equation}
with $\varphi(x)=\varphi_0(x)=-\log x+x-1$. The explicit form of
$m(\theta,t)$ in this case is
\begin{eqnarray}\label{m x the pr klm}
    x \mapsto m(x,\theta,t) & = & t_{0}-\left(t^T\overline{g}(x,\theta)\right)
    \frac{1}{1-t^T\overline{g}(x,\theta)}+\log\left(1-t^T\overline{g}(x,\theta)\right)
    +\frac{1}{1-t^T\overline{g}(x,\theta)}-1.\nonumber\\
    & = & t_0+\log\left(1-t^T\overline{g}(x,\theta)\right).
\end{eqnarray}
For fixed $\theta\in\Theta$, the sup in (\ref{klm M theta}), which
we have noted $\widehat{c_\theta}$, satisfies the following system
\begin{equation}\label{systel l+1 eq l+1 var}
\left\{
\begin{array}{ccl}
  \int \frac{1}{1-c_0-\sum_{j=1}^lc_jg_j(x,\theta)}~dP_n(x) & = & 1 \\
  \int \frac{g_j(x,\theta)}{1-c_0-\sum_{j=1}^lc_jg_j(x,\theta)}~dP_n(x)
  & = & 0,~~\text{ for all } j=1,...,l \\
\end{array}
\right.
\end{equation}
a system of $(l+1)$ equations and $(l+1)$ variables. The
projection $\widehat{Q^*_\theta}$ is then obtained using
Proposition \ref{existence de Q star n} part (ii). We have for all
$i\in\left\{1,\ldots,n\right\}$
\begin{equation*}
\frac{1}{\widehat{Q^*_\theta}(X_i)}=n\left(1-c_0-\sum_{j=1}^nc_jg_j(X_i,\theta)\right)
\end{equation*}
which, multiplying by $\widehat{Q^*_\theta}(X_i)$ and summing upon
$i$ yields $c_0=0$. Therefore the system (\ref{systel l+1 eq l+1
var}) reduces to the system (3.3) in \cite{Qin-Lawless1994}
replacing $c_1,\ldots,c_l$ by $-t_1,\ldots,-t_l$. Simplify
 (\ref{m x the pr klm}) plugging $t_0=0$. Notice that
$2n\widehat{KL_m}(\mathcal{M}_\theta,P_0)=l_E(\theta_0)$ in the
notation of \cite{Qin-Lawless1994}, and that the function of
$t=(0,-\tau_1,\ldots,-\tau_l)$ defined by
\begin{equation*}
    t \mapsto P_nm(\theta,t)
\end{equation*}
coincide with the function
\begin{equation*}
    \tau \rightarrow P_n\log\left(1+\tau^Tg(.,\theta)\right)
\end{equation*}
used in \cite{Qin-Lawless1994}. The interest in formula (\ref{klm
M theta}) lays in the obtention of the limit distributions of
$2n\widehat{KL_m}(\mathcal{M}_\theta,P_0)$ under $\mathcal{H}_1$.
By Theorem \ref{lois limite a theta fixe}, we have
\begin{equation*}
    \sqrt{n}\left(\widehat{KL_m}(\mathcal{M}_\theta,P_0)-
    KL_m(\mathcal{M}_\theta,P_0)\right)
\end{equation*}
converges to a normal distribution variable, which proves
consistency of the test; this results cannot be obtained by
the  \cite{Qin-Lawless1994}'s approach.\\

\noindent The choice of $\varphi$ depends on some a priori
knowledge on $\theta_0$. Hopefully, some divergences do not have
such an inconvenient. We now clarify this point. For fixed
$\theta$ in $\Theta$, let $\mathcal{M}_\theta^{(n)}$ and
$\mathcal{D}_\phi^{(n)}$ be defined respectively as in (\ref{m
theta n}) and in (\ref{D phi n}).
 Assume that
$\mathcal{M}_\theta^{(n)}\cap\mathcal{D}_\phi^{(n)}$ is not void.
Then $P_n$ has a projection $\widehat{Q_\theta^*}$ on
$\mathcal{M}_\theta^{(n)}$ and $\phi(\widehat{Q_\theta^*},P_n)$ is
finite. The estimation of $\theta_0$ is achieved minimizing
$\widehat{\phi}(\mathcal{M}_\theta,P_0)$ on the sets
\begin{equation*}
    \Theta_n^\phi:=\left\{\theta\in\Theta \text{ such that }
     \mathcal{M}_\theta^{(n)}\cap\mathcal{D}_\phi^{(n)}~\text{ is not void}
     \right\}.
\end{equation*}
Clearly the description of $\Theta_n^\phi$ depends on the
divergence $\phi$. Consider the following example, with $n=2$,
$X=(X_1,X_2)$ and $g(x,\theta)=x-\theta$. Then
\begin{equation*}
    \mathcal{M}_\theta=\left\{(q_1,q_2)^T~\text{ such that  }~ q_1+q_2=1
    \text{ and }q_1(X_1-\theta)+q_2(X_2-\theta)=0\right\}
\end{equation*}
and
\begin{equation*}
    \mathcal{D}^{(2)}_\phi=\left\{(q_1,q_2) ~ \text{ such that }~ \frac{1}{2}
    \sum_{i=1}^2\varphi(2q_i)<\infty\right\}.
\end{equation*}
When $\phi=KL_m$, then
$\mathcal{D}^{(2)}_{KL_m}=\mathbb{R}_+^*\times\mathbb{R}_+^*$. So,
according to the value of $\theta$,
$\mathcal{M}_\theta^{(n)}\cap\mathcal{D}_{KL_m}^{(n)}$ may be void
and therefore $\Theta_n^{KL_m}$ has a complex structure. At the
opposite, for example when $\phi=\chi^2$, then
$\mathcal{D}^{(2)}_{\chi^2}=\mathbb{R}^2$. Hence
$\mathcal{M}_\theta^{(n)}\cap\mathcal{D}_\phi^{(n)}=
\mathcal{M}_\theta^{(n)}$ which is not void for all $\theta$ and
hence $\Theta_n^{\chi^2}=\Theta$.\\

\noindent On the other hand, we have for any $\phi$-divergence
\begin{eqnarray*}
\widehat{\theta}_\phi & := &
\arg\inf_{\theta\in\Theta}\inf_{Q\in\mathcal{M}_\theta^{(n)}}\widehat{\phi}
(Q,P_0)\\
 & = & \arg\inf_{\theta\in\Theta}\inf_{Q\in\mathcal{M}_\theta^{(n)}
 \cap \mathcal{D}_\phi^{(n)}}\widehat{\phi}(Q,P_0).
\end{eqnarray*}
When $\mathcal{D}_\phi^{(n)}\neq\mathbb{R}^n$, the infimum in
$\theta$ above should be taken upon $\Theta_n^\phi$ which might be
quite cumbersome.  \cite{Owen2001}  indeed mentions such a
difficulty.\\

\noindent In relation to this problem, \cite{Qin-Lawless1994}
bring some asymptotic arguments in the case of the empirical
likelihood. They show that there exists a sequence of
neighborhoods
\begin{equation*}
    V_n(\theta_0):=\left\{\theta ~\text{ such that }~\|\theta-\theta_0\|\leq
    n^{-1/3}\right\}
\end{equation*}
on which, with probability one as $n$ tends to infinity,
$L_n(\theta)$ has a maximum. This turns out, in the context of
$\phi$-divergences, to write that the mapping
\begin{equation*}
    \theta \mapsto \inf_{Q\in\mathcal{M}_\theta^{(n)}}KL_m(Q,P_n)
\end{equation*}
has a minimum when $\theta$ belongs to $V_n(\theta_0)$. This
interesting result does not solve the problem for fixed $n$, as
$\theta_0$ is unknown. For such problem, the use of
$\phi$-divergences, satisfying
$\mathcal{D}_\phi^{(n)}=\mathbb{R}^n$ (for example
$\chi^2$-divergence), might give information about $\theta_0$ and
localizes it through
$\phi$-divergence confidence regions ($CR_\phi$'s).\\

\noindent The choice of the divergence $\phi$ also depends upon
some knowledge on the support of the unknown p.m. $P_0$. When
$P_0$ has a projection on $\mathcal{M}$ with same support as
$P_0$, Proposition \ref{existence de Q star}
 yields its description and its explicit calculation.
A necessary condition for this is that $\mathcal{C}_\theta$, as
defined in (\ref{C theta}), has non void interior in
$\mathbb{R}^{(l+1)}$. Consider the case of the empirical
likelihood, that is when $\varphi(x)=-\log x+x-1$; then $\text{ Im
}\varphi'=]-\infty,1[$. Consider $g(x,\theta)=x-\theta$, i.e., a
constraint on the mean. Assume that the support of $P_0$ is
unbounded. Then
\begin{equation*}
    \mathcal{C}_\theta=\left\{t\in\mathbb{R}^2
    \text{ such that  for all } x ~(P_0 - a.s.)~,~t_0+t_1(x-\theta)\in ]-\infty,1[
    \right\}.
\end{equation*}
Therefore, $t_1=0$ and $\mathcal{C_\theta}=]-\infty,1[\times\{0\}$
which implies that the interior of $\mathcal{C}_\theta$ is void.
This results indicates that the support of $Q^*$ is not the same
as the support of $P_0$. Hence in this case we cannot use the dual
representation of $KL_m(\mathcal{M}_\theta,P_0)$. The arguments
used in Section 3 for the obtention of limiting distributions
cannot be used, if the support of $P_0$ is unbounded, in order to
obtain the limiting distribution of the estimates
$\widehat{KL_m}(\mathcal{M}_\theta,P_0)$ under $\mathcal{H}_1$
(i.e., when $P_0$ does not belong to $\mathcal{M}_\theta$). We
thus cannot conclude in this case that the tests pertaining to
$\theta_0$ are consistent.

\section{Robustness and Efficiency of ME$\phi$D estimates and Simulation Results}

\noindent  \cite{Lindsay1994} introduced a general instrument for
the study of the asymptotic properties of parametric estimates by
minimum $\phi$-divergences, called Residual Adjustment Function
(RAF). We first recall its definition. Let
$\left\{P_\theta~;~\theta\in\Theta\right\}$ be some parametric
model defined on a finite set $\mathcal{X}$. Let $X_1,\ldots,X_n$
a sample with distribution $P_{\theta_0}$. A minimum
$\phi$-divergence estimate (M$\phi$DE) (called also minimum
disparity estimator) of $\theta_0$ is given by
\begin{equation}\label{minimum phi div estim param}
    \widetilde{\theta}_\phi:=\arg\inf_{\theta\in\Theta}\sum_{x\in\mathcal{X}}
    \varphi\left(\frac{P_\theta(x)}{P_n(x)}\right)P_n(x),
\end{equation}
where $P_n(x)$ is the proportion of the sample point that take
value $x$. When the parametric model
$\left\{P_\theta~:~\theta\in\Theta\right\}$ is regular, then
$\widetilde{\theta}_\phi$ is solution of the equation
\begin{equation}\label{sum var phi prime}
\sum_{x\in\mathcal{X}}
    \varphi'\left(\frac{P_\theta (x)}{P_n (x)}\right)\dot{P}_\theta
    (x)=0,
\end{equation}
which can be written as
\begin{equation}\label{A phi de x}
\sum_{x\in\mathcal{X}} A_\varphi(\delta(x))\dot{P}_\theta(x)=0.
\end{equation}
 In this display, $A_\varphi(u):=\varphi'\left(\frac{1}{u+1}\right)$
 depends only upon the divergence function $\varphi$ and
\begin{equation*}
\delta(x):=\frac{P_n(x)}{P_\theta(x)}-1
\end{equation*}
is the ``Pearson Residual'' at $x$ which belongs to
$]-1,+\infty[$. The function $A_\varphi(.)$ is the RAF.\\

\noindent The points $x$ for which $\delta(x)$ is close to $-1$
are called ``inliers'', whereas points $x$ such that $\delta(x)$
is large are called ``outliers''. Efficiency properties are linked
with the behavior of $A_\varphi(.)$ in the neighborhood of $0$
(see \cite{Lindsay1994} Proposition 3 and \cite{Basu-Lindsay1994})
: the smaller the value of $\left|A_\varphi''(0)\right|$, the more
second efficient the estimate $\widetilde{\theta}_\phi$ in the
sense of \cite{Rao1961}.\\

\noindent It is easy to verify that the RAF's of the power
divergences $\phi_\gamma$, defined by the divergence functions in
(\ref{gamma convex functions}), have the form
\begin{equation}\label{A gamma de delta}
A_\gamma(\delta) = \frac{\left(\delta
+1\right)^{1-\gamma}-1}{(\gamma -1)}.
\end{equation}
In particular, the M$\phi_\gamma$DE of (\ref{sum var phi prime})
with the RAF in (\ref{A gamma de delta}) corresponds to the
maximum likelihood when $\gamma=0$, minimum Hellinger distance
when $\gamma=0.5$, minimum $\chi^2$ divergence when $\gamma = 2$,
minimum modified $\chi^2$ divergence when $\gamma = -1$ and
minimum $KL$ divergence
when $\gamma=1$.\\

\noindent From (\ref{A gamma de delta}), we see that
$A_\gamma''(0)=\gamma$. Hence for the maximum likelihood estimate,
we have $\left|A_\gamma''(0)\right|=\left|A_0''(0)\right|=0$ which
is the smallest value of $\left|A_\gamma''(0)\right|$,
$\gamma\in\mathbb{R}$. Therefore, according to Proposition 3 in
\cite{Lindsay1994}, the maximum likelihood estimate is the most
second-order efficient estimate (in the sense of \cite{Rao1961})
among all minimum power divergences estimates.\\

\noindent Robustness features of $\widetilde{\theta}_\phi$ against
inliers and outliers are related to the variations of
$A_\varphi(u)$ or $\varphi(x)$ when $u$ or $x$ close to $-1$ and
$+\infty$, respectively as seen through the following heuristic
arguments. Let $\phi_1$ and $\phi_2$ two divergences associated to
the functions $\varphi_1$ and $\varphi_2$. If
\begin{equation*}
    \lim_{x \downarrow
    0}\frac{\varphi_1(x)}{\varphi_2(x)}=+\infty,
\end{equation*}
then the estimating equation (\ref{sum var phi prime})
corresponding to $\varphi_1$ in not as stable as that
corresponding to $\varphi_2$, and hence the ME$\phi_2$DE is more
robust than ME$\phi_1$DE against outliers. If
\begin{equation*}
    \lim_{x \uparrow
    +\infty}\frac{\varphi_1(x)}{\varphi_2(x)}=+\infty,
\end{equation*}
then the estimating equation (\ref{sum var phi prime})
corresponding to $\varphi_1$ is not as stable as that
corresponding to $\varphi_2$, and hence the ME$\phi_2$DE is more
robust than ME$\phi_1$DE against inliers.\\

\noindent In all cases, the divergence associated to the
divergence function having the smallest variations on its domain
leads to the most
robust estimate against both outliers and inliers.\\

\noindent It is shown also in \cite{JimenezShao2001} that no
minimum power divergence estimate (including the maximum
likelihood one) is better than the minimum Hellinger divergence in
terms of both
second-order efficiency and robustness.\\

\noindent In the examples below, we compare by simulations the
efficiency and robustness properties of some ME$\phi$DE's for some
models satisfying linear constraints. We will see that the minimum
empirical Hellinger divergence estimate represents a suitable
compromise between efficiency and robustness.  A theoretical study
of efficiency and robustness properties of ME$\phi$DE's is
necessary and should envolve second-order efficiency versus
robustness since all ME$\phi$DE's are all equally first-order
efficient (see Remark  \ref{tous les estim sont equiv} and Theorem
\ref{loi limite de thetan et cn}).

\subsection*{Numerical  Results}

We consider for illustration the same model as in
\cite{Qin-Lawless1994} Section 5 Example 1. The model
$\mathcal{M}_\theta$ (see \ref{modele signees finies}) here is the
set of all signed finite measures $Q$ satisfying
\begin{equation}\label{modele simulations}
    \int dQ = 1 ~\text{ and } ~ \int g(x,\theta)~dQ(x)=0,
\end{equation}
with $g(x,\theta)=\left((x-\theta),(x^2-2\theta^2-1)\right)^T$ and
$\theta$, the parameter of interest, belongs to $\mathbb{R}$.\\

\noindent In Examples 1.a and 1.b below, we compare the efficiency
property of various estimates : we generate $1000$ pseudorandom
samples of sizes 25, 50, 75 and 100 from a normal distribution
with mean $\theta_0$ and variance $\theta_0^2+1$ (i.e.,
$P_0=\mathcal{N}(\theta_0,\theta_0^2+1)$) for two values of
$\theta_0$ : $\theta_0=0$ in Example 1.a and $\theta_0=1$ in
Example 1.b. Note that $P_0$ satisfies
(\ref{modele simulations}).\\

\noindent For each sample, we consider various estimates of
$\theta_0$ : the sample mean estimate (SME), the parametric ML
estimate (MLE) based on the normal distribution
$\mathcal{N}(\theta,\theta^2+1)$ and ME$\phi$D estimates
$\widehat{\theta}_\phi$ associated to the divergences :
$\phi=\chi^2_m$, $H$, $KL$, $\chi^2$ and $KL_m$-divergence (which
coincides with the MEL one, i.e., ME$KL_m$E=MELE).\\

\noindent For all divergence $\phi$ considered, in order to
calculate the ME$\phi$DE $\widehat{\theta}_\phi$,  we first
calculate $\widehat{\phi}(\mathcal{M}_\theta,P_0)$ for all given
$\theta$ (using the representation (\ref{estim de phiMtheta
dual})) by Newton's method, and then minimize it to obtain
$\widehat{\theta}_\phi$.\\

\noindent The results of Theorem \ref{loi limite de thetan et cn}
show that for all $\phi$-divergence
\begin{equation*}
    \sqrt{n}\left(\widehat{\theta}_\phi-\theta_0\right)
    \rightarrow \mathcal{N}(0,V)
\end{equation*}
where $V$ is independent of the divergence $\phi$; it is given in
Theorem \ref{loi limite de thetan et cn}. For the present model,
following \cite{Qin-Lawless1994}, $V$ writes
\begin{equation}\label{relation entre les var}
V=Var(X)-\triangle^{-1}\left[m'(\theta_0)Var(X)+\theta_0m(\theta_0)-E(X^3)\right])^2
\end{equation}
where
$\triangle=E\left[m'(\theta_0)(X-\theta_0)+m(\theta_0)-X^2\right]^2$
and $m(\theta):=2\theta^2+1$. Thus $V\leq Var(X)$ which is the
variance of $\sqrt{n}\left(\overline{X}_n-\theta_0\right)$ with
$\overline{X}_n:=\frac{1}{n}\sum_{i=1}^{n}X_i$, the sample mean
estimate (SME) of $\theta_0$. So, EM$\phi$D estimates are all
asymptotically at least as efficient as $\overline{X}_n$.

\subsection{Example 1.a}

In this example the true value of the parameter is $\theta_0=0$.

\vskip 0.5cm

\begin{table}[h]
\begin{tabular}{ || c  ||c c|c c|c c|c c||}
 \hline \hline   & \multicolumn{2}{c|}{ME$\chi^2_m$DE} &
 \multicolumn{2}{c|}{ME$KL_m$DE=MELE} &
 \multicolumn{2}{c|}{ME$H$DE} & \multicolumn{2}{c||}{ME$KL$DE} \\
 $n$  &  mean   &  var    &  mean   &  var   &  mean    &  var   &  mean   &  var   \\ \hline
 25   & 0.0089  & 0.0314  & 0.0086  & 0.0315 & 0.0084   & 0.0315 & 0.0082  & 0.0314  \\
 50   & -0.0116 & 0.0209  & -0.0118 & 0.0210 & -0.0119  & 0.0210 & -0.0120 & 0.0210  \\
 75   & -0.0025 & 0.0171  & -0.0024 & 0.0170 & -0.0023    & 0.0170 & -0.0022 & 0.0169 \\
 100  & -0.0172 & 0.0112  & -0.0174 & 0.0111 & -0.0174  & 0.0111 & -0.0175 & 0.0112 \\
\hline \hline & \multicolumn{2}{c|}{ME$\chi^2$DE} &
 \multicolumn{2}{c|}{PMLE} &
 \multicolumn{2}{c|}{SME}  & \multicolumn{2}{c||}{ } \\
 $n$  &  mean   &  var    &  mean   &  var   &  mean    &  var   &   & \\
 \hline
 25   & 0.0077  & 0.0313  & 0.0026  & 0.0318 & 0.0081   & 0.0394 &   &   \\
 50   & -0.0125 & 0.0212  & -0.0063 & 0.0196 & -0.0040  & 0.0200 &   &   \\
 75   & -0.0019 & 0.0167  & -0.0011 & 0.0170 & 0.0013   & 0.0164 &   &   \\
 100  & -0.0177 & 0.0112  & -0.0158 & 0.0108 & -0.0149  & 0.0102 &   &   \\
 \hline \hline
\end{tabular}
\vskip 0.5cm \caption{Estimated mean and variance of the estimates
of $\theta_0$ in Example 1.a.} \label{Efficacitetheta0-0}
\end{table}

\noindent We can see from Table \ref{Efficacitetheta0-0} that all
the estimates converge in a satisfactory way. The estimated
variances are almost the same for all estimates. This is not
surprising since the limit variance  of all estimates in this
Example (when $\theta_0=0$) is close to $V(X)$.

\subsection{Example 1.b}

In this example the true value of the parameter is $\theta_0=1$.

\vskip 0.5cm

\begin{table}[h]
\begin{tabular}{ || c  ||c c|c c|c c|c c||}
 \hline \hline   & \multicolumn{2}{c|}{ME$\chi^2_m$DE} &
 \multicolumn{2}{c|}{ME$KL_m$DE=MELE} &
 \multicolumn{2}{c|}{ME$H$DE} & \multicolumn{2}{c||}{ME$KL$DE} \\
 $n$  &  mean   &  var    &  mean   &  var   &  mean    &  var   &  mean   &  var   \\ \hline
 25   & 0.9394  & 0.0310  & 0.9387  & 0.0312 & 0.9385   & 0.0313 & 0.9378  & 0.0316  \\
 50   & 0.9994  & 0.0186  & 0.9967  & 0.0186 & 0.9954   & 0.0186 & 0.9941  & 0.0187  \\
 75   & 1.0009  & 0.0156  & 0.9988  & 0.0154 & 0.9975   & 0.0154 & 0.9966  & 0.0153 \\
 100  & 0.9984  & 0.0113  & 0.9959  & 0.0112 & 0.9945   & 0.0112 & 0.99315 & 0.0112 \\
\hline \hline & \multicolumn{2}{c|}{ME$\chi^2$DE} &
 \multicolumn{2}{c|}{PMLE} &
 \multicolumn{2}{c|}{SME}  & \multicolumn{2}{c||}{ } \\
 $n$  &  mean   &  var    &  mean   &  var   &  mean    &  var   &   & \\
 \hline
 25   & 0.9350  & 0.0322  & 0.9540  & 0.0325 & 1.0033   & 0.0810 &   &   \\
 50   & 0.9909  & 0.0190  & 1.0036  & 0.0174 & 1.0021   & 0.0407 &   &   \\
 75   & 0.9940  & 0.0152  & 1.0003  & 0.0149 & 0.9912   & 0.0288 &   &   \\
 100  & 0.9900  & 0.0113  & 0.9970  & 0.0107 & 0.9851   & 0.0262 &   &   \\
 \hline \hline
\end{tabular}
\vskip 0.5cm \caption{Estimated mean and variance of the estimates
of $\theta_0$ in Example 1.b.} \label{Efficacitetheta0-1}
\end{table}

\begin{figure}[h]
  \centerline{
    \includegraphics[width=.55\textwidth]{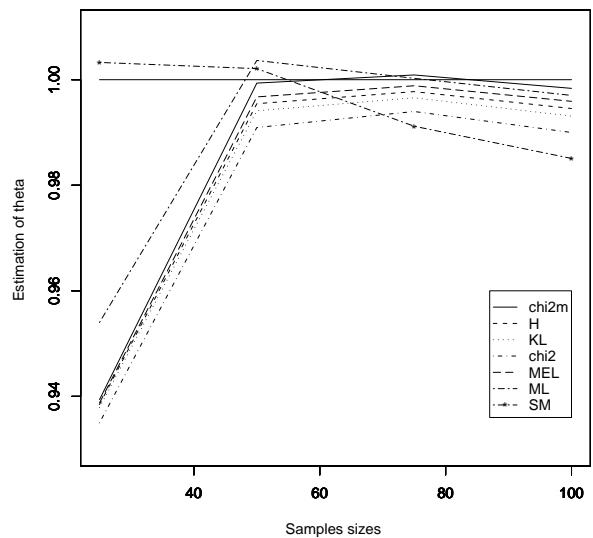}}
  \caption{Estimated mean of the estimates of $\theta_0$ in Example 1.b.}
  \label{mean dans l exemple 1.b}
\end{figure}

\noindent We can see from Table \ref{Efficacitetheta0-1} and
Figure \ref{mean dans l exemple 1.b} that  the estimated bias of
E$\phi$DE's are all smaller than the SME one for moderate and
large sample sizes. Furthermore, from Figure \ref{var dans l
exemple 1.b}, we observe that the estimated variances of
E$\phi$DE's are all less than the SME one. They lie between that
of the sample mean and  that of the parametric maximum likelihood
estimate. We observe also that the estimated variances of the MELE
and ME$H$DE are equal and are the smallest among the variances of
all ME$\phi$DE's considered. It should be emphasized that even for
small sample sizes, the MSE of the SM is larger than any of
ME$\phi$DE's.
\begin{figure}[h]
  \centerline{
    \includegraphics[width=.55\textwidth]{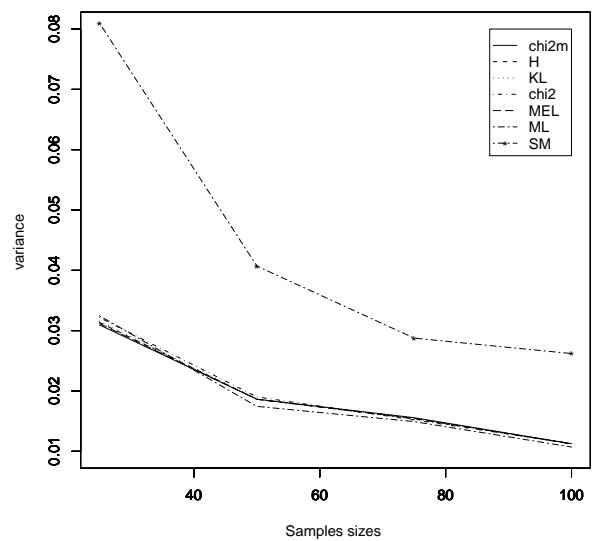}}
  \caption{Estimated variance of the estimates of $\theta_0$ in Example 1.b.}
  \label{var dans l exemple 1.b}
\end{figure}

\noindent In Examples 2.a and 2.b below, we compare robustness
property of the estimates considered above for contaminated data :
we consider the same model $\mathcal{M}_\theta$ as in (\ref{modele
simulations}).

\subsection{Example 2.a}

In this Example, we generate $1000$ pseudo-random samples of sizes
25, 50, 75 and 100 from  a distribution
\begin{equation*}
  \widetilde{P_0} =   (1-\epsilon)P_0+\epsilon\delta_5
\end{equation*}
where $P_0=\mathcal{N}(\theta_0,\theta_0^2+1)$, $\epsilon=0.15$
and $\theta_0 = 2$.  We consider the same estimates as in the
above examples.

\vskip 0.5 cm

\begin{table}[h]
\begin{tabular}{ || c  ||c c|c c|c c|c c||}
 \hline \hline   & \multicolumn{2}{c|}{ME$\chi^2_m$DE} &
 \multicolumn{2}{c|}{ME$KL_m$DE=MELE} &
 \multicolumn{2}{c|}{ME$H$DE} & \multicolumn{2}{c||}{ME$KL$DE} \\
 $n$  &  mean   &  var    &  mean   &  var   &  mean    &  var   &  mean   &  var   \\ \hline
 25   & 2.1609  & 0.0654  & 2.1513  & 0.0653 & 2.1453   & 0.0653 & 2.1396  & 0.0652  \\
 50   & 2.2087  & 0.0303  & 2.1975  & 0.0304 & 2.1912   & 0.0307 & 2.1848  & 0.0309  \\
 75   & 2.2218  & 0.0214  & 2.2106  & 0.0213 & 2.2046   & 0.0213 & 2.1987  & 0.0215 \\
 100  & 2.2283  & 0.0151  & 2.2169  & 0.0149 & 2.2110   & 0.0148 & 2.2052  & 0.0149 \\
\hline \hline & \multicolumn{2}{c|}{ME$\chi^2$DE} &
 \multicolumn{2}{c|}{PMLE} &
 \multicolumn{2}{c|}{SME}  & \multicolumn{2}{c||}{ } \\
 $n$  &  mean   &  var    &  mean   &  var   &  mean    &  var   &   & \\
 \hline
 25   & 2.1278  & 0.0646  & 2.2088  & 0.0581 & 2.4265   & 0.2178 &   &   \\
 50   & 2.1729  & 0.0316  & 2.2296  & 0.0280 & 2.4535   & 0.1076 &   &   \\
 75   & 2.1877  & 0.0219  & 2.2337  & 0.0197 & 2.4545   & 0.0721 &   &   \\
 100  & 2.1947  & 0.0151  & 2.2352  & 0.0139 & 2.4572   & 0.0543 &   &   \\
 \hline \hline
\end{tabular}
\vskip 0.5cm \caption{Estimated mean and variance of the estimates
of $\theta_0$ in Example 2.a.} \label{Robustessetheta0-2-Outlier5}
\end{table}

\noindent In this Example, we can see from Table
\ref{Robustessetheta0-2-Outlier5} and Figure \ref{mean dans l
exemple 2.a} that the ME$\chi^2$D estimate is the most robust and
ME$\chi^2_m$ estimate is the least robust. We observe also that
the MELE which is the ME$KL_m$DE is less robust than the ME$KL$DE
and that the ME$H$D estimate is more robust than MEL one.

\begin{figure}[h]
  \centerline{
    \includegraphics[width=.55\textwidth]{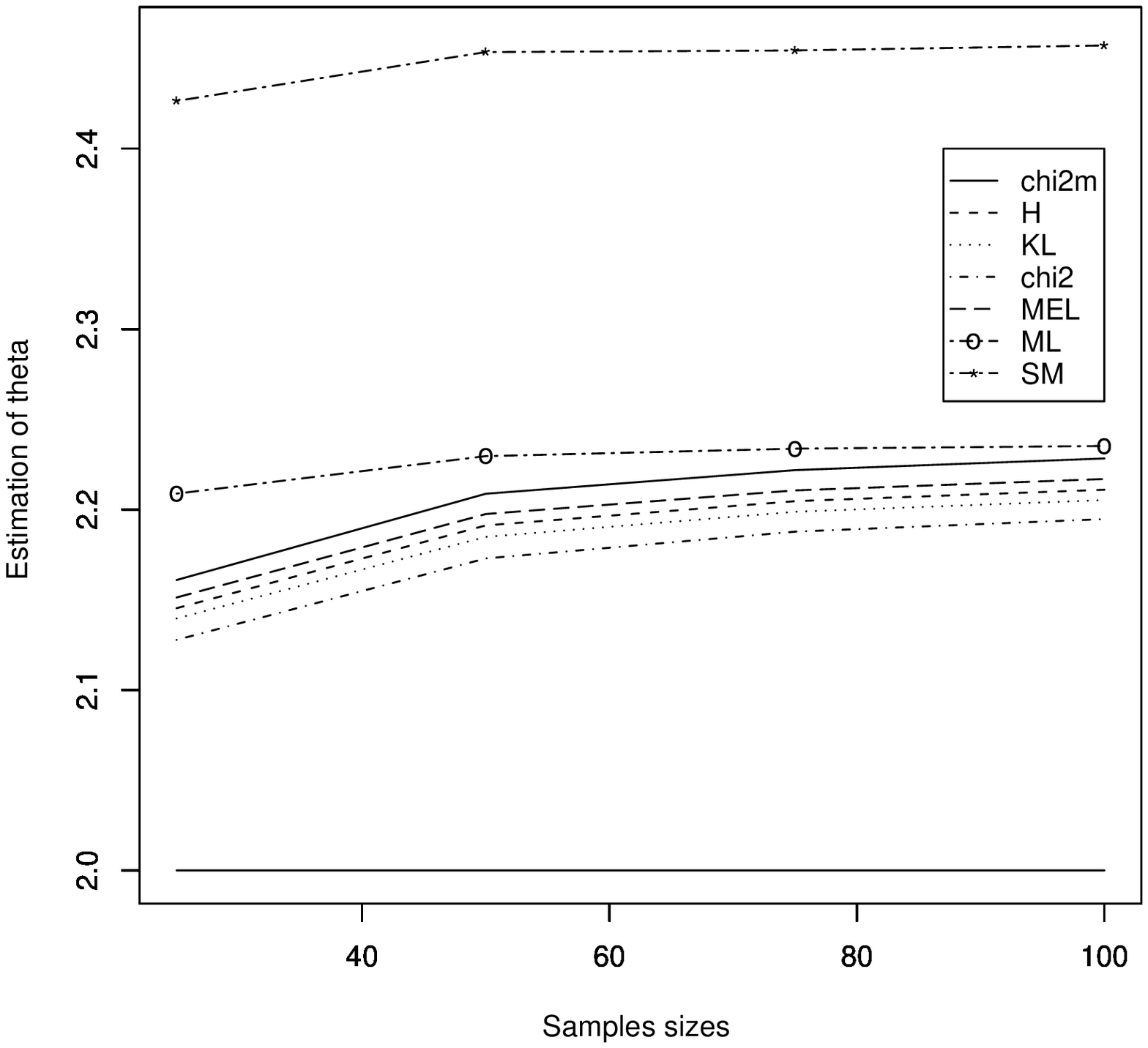}}
  \caption{Estimated mean  of the estimates of $\theta_0$ in Example 2.a.}
  \label{mean dans l exemple 2.a}
\end{figure}

\subsection{Example 2.b}

In this Example, we generate $1000$ pseudo-random samples of sizes
50, 100, 150 and 200  from a distribution
$P_0=\mathcal{N}(\theta_0,\theta_0^2+1)$ with $\theta_0 = 2$ and
we cancel the observations  in the interval $[4,5]$ . We consider
the same estimates as in the above examples.

\vskip 0.5cm

\begin{table}[h]
\begin{tabular}{ || c  ||c c|c c|c c|c c||}
 \hline \hline   & \multicolumn{2}{c|}{ME$\chi^2_m$DE} &
 \multicolumn{2}{c|}{ME$KL_m$DE=MELE} &
 \multicolumn{2}{c|}{ME$H$DE} & \multicolumn{2}{c||}{ME$KL$DE} \\
 $n$  &  mean   &  var    &  mean   &  var   &  mean    &  var   &  mean   &  var   \\ \hline
 50   & 1.9917  & 0.0451  & 1.9784  & 0.0431 & 1.9721   & 0.0426 & 1.9659  & 0.0423  \\
 100  & 1.9962  & 0.0362  & 1.9844  & 0.0346 & 1.9787   & 0.0341 & 1.9729  & 0.0336  \\
 150  & 2.0011  & 0.0150  & 1.9903  & 0.0142 & 1.9849   & 0.0139 & 1.9795  & 0.0137 \\
 200  & 1.9602  & 0.0162  & 1.9516  & 0.0158 & 1.9473   & 0.0157 & 1.9430  & 0.0156 \\
\hline \hline & \multicolumn{2}{c|}{ME$\chi^2$DE} &
 \multicolumn{2}{c|}{PMLE} &
 \multicolumn{2}{c|}{SME}  & \multicolumn{2}{c||}{ } \\
 $n$  &  mean   &  var    &  mean   &  var   &  mean    &  var   &   & \\
 \hline
 50   & 1.9522  & 0.0428  & 1.9705  & 0.0358 & 1.7750   & 0.1039 &   &   \\
 100  & 1.9590  & 0.0329  & 1.9687  & 0.0298 & 1.7365   & 0.0576 &   &   \\
 150  & 1.9671  & 0.0135  & 1.9781  & 0.0121 & 1.7456   & 0.0283 &   &   \\
 200  & 1.9325  & 0.0155  & 1.9420  & 0.0146 & 1.7247   & 0.0317 &   &   \\
 \hline \hline
\end{tabular}
\vskip 0.5cm \caption{Estimated mean and variance of the estimates
of $\theta_0$ in Example 2.b.}
\label{Robustessetheta0-2-Inlier4-5}
\end{table}

\noindent In this example, in contrast with Example 2.b,  we
observe that the ME$\chi^2_m$DE is the most robust,  ME$\chi^2$DE
is the least robust and ME$KL$DE is less robust than ME$KL_m$DE
(=MELE). Generally, if a ME$\phi$DE is more robust than its
adjoint\footnote{For all divergence $\phi$ associated to a convex
function $\varphi$, its adjoint, noted $\phi^\sim$, is the
divergence associated to the convex function, noted
$\varphi^\sim$, defined by : $\varphi^\sim (x)=x\varphi(1/x)$, for
all $x$.} (i.e., ME$\phi^\sim$DE) against ``outliers'', then it is
less robust then its adjoint against ``inliers'' (see Examples 2.a
and 2.b). The Hellinger divergence has not this disadvantage since
it is self-adjoint (i.e., $H=H^\sim$).

\begin{figure}[h]
  \centerline{
    \includegraphics[width=.55\textwidth]{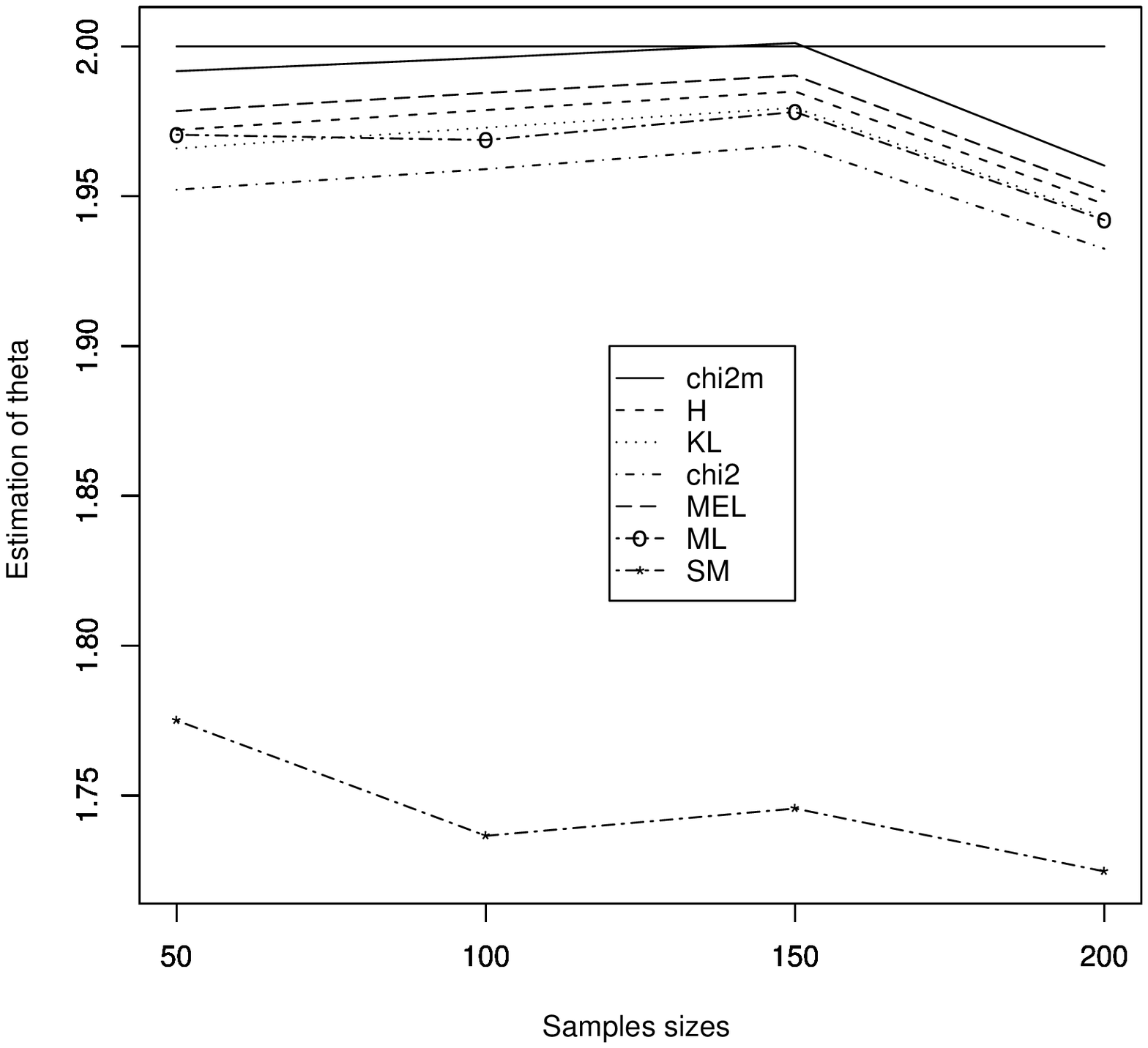}}
  \caption{Estimated mean  of the estimates of $\theta_0$ in Example 2.b.}
  \label{mean dans l exemple 2.b}
\end{figure}

\section{Proofs}

\subsection{Proof of  Proposition \ref{existence de Q star n}}
Proof of part (i). The function
\begin{equation*}
\left(Q(X_1),\ldots,Q(X_n)\right)^T\in\mathbb{R}^n\mapsto
\frac{1}{n}\sum_{i=1}^{n}\varphi\left(nQ(X_i)\right)
\end{equation*}
is continuous and nonnegative on $\mathcal{D}_\phi^{(n)}$.
Furthermore, the set $\mathcal{M}_\theta^{(n)}$ is closed in
$\mathbb{R}^n$. Hence, by condition (\ref{condition d'existence de
Q star n}), the infimum of the function
\begin{equation*}
\left(Q(X_1),\ldots,Q(X_n)\right)^T\in\mathbb{R}^n\mapsto
\frac{1}{n}\sum_{i=1}^{n}\varphi\left(nQ(X_i)\right)
\end{equation*}
on the set $\mathcal{D}_\phi^{(n)}\cap \mathcal{M}_\theta^{(n)}$
exists as an interior point of $\mathcal{D}_\phi^{(n)}$. Since the
above function is strictly convex and the set
$\mathcal{D}_\phi^{(n)}\cap \mathcal{M}_\theta^{(n)}$ is convex,
then this infimum is unique. It is noted $\widehat{Q_\theta^*}$.
This
concludes the proof of part (i).\\
Proof of part (ii). Since
$\left(Q(X_1),\ldots,Q(X_n)\right)^T\in\mathbb{R}^n\mapsto
\frac{1}{n}\sum_{i=1}^{n}\varphi\left(nQ(X_i)\right)$ is
$\mathcal{C}^1$ on the interior of $\mathcal{D}_\phi^{(n)}$, and
since $\widehat{Q_\theta^*}$ is in the interior of
$\mathcal{D}_\phi^{(n)}$, we can use the Lagrange method. This
yields the explicit form (\ref{caracterisation de Q star n}) of
the projection $\widehat{Q_\theta^*}$ in which $\widehat{c}_0$ is
the Lagrange multiplier associated to the constraint
$\sum_{i=1}^{n}Q(X_{i})=1$ and $\widehat{c}_j$ to the constraint
$\sum_{i=1}^{l}Q(X_{i})g_{j}(X_{i},\theta )=0$, for all
$j=1,\ldots,l$. This concludes the proof of  Proposition
\ref{existence de Q star n}.

\subsection{Proof of Proposition \ref{consitence de l estimateur de
ctheta}}

Define the estimates
\begin{equation*}
    \widetilde{c_\theta}=\arg\inf_{t\in T_\theta} P_nm(\theta,t)
    ~\text{ and }~
    \widetilde{\phi}\left(\mathcal{M}_\theta,P_0\right)=
    \sup_{t\in T_\theta} P_nm(\theta,t).
\end{equation*}
 By condition
(C.2), for all $n$ sufficiently large, we have
\begin{equation*}
    \widehat{c_\theta}=\widetilde{c_\theta} ~\text{ and }~
        \widehat{\phi}\left(\mathcal{M}_\theta,P_0\right)=
            \widetilde{\phi}\left(\mathcal{M}_\theta,P_0\right).
\end{equation*}
We prove that
$\widetilde{\phi}\left(\mathcal{M}_\theta,P_0\right)$ and
$\widetilde{c_\theta}$  converge to
$\phi\left(\mathcal{M}_\theta,P_0\right)$ and $c_\theta$
respectively. Since $c_\theta$ is isolated, then consistency of
$\widetilde{c_\theta}$ holds as a consequence of Theorem 5.7 in
\cite{vanderVaart1998}. For the estimate
$\widetilde{\phi}\left(\mathcal{M}_\theta,P_0\right)$, we have
\begin{equation*}
\left|\widetilde{\phi}\left(\mathcal{M}_\theta,P_0\right)-
\phi\left(\mathcal{M}_\theta,P_0\right)\right|=\left|
P_nm(\theta,\widetilde{c_\theta})-P_0m(\theta,c_\theta)\right|:=|A|,
\end{equation*}
which implies
\begin{equation*}
    P_nm(\theta,c_\theta)-P_0m(\theta,c_\theta)< A <
    P_n
    m(\theta,\widetilde{c_\theta})-P_0m(\theta,\widetilde{c_\theta}).
\end{equation*}
Both the RHS and the LHS terms in the above display go to $0$,
under condition (C.1). This implies that $A$ tends to $0$, which
concludes the proof of Proposition \ref{consitence de l estimateur
de ctheta}.

\subsection{Proof of Theorem \ref{lois limite a theta fixe}}.
Proof of part (1). Some calculus yield
\begin{eqnarray}\label{condition de 1o}
    P_0m'(\theta,c_\theta) & = & P_0\left(1-\overleftarrow{\varphi'}
    \left(c_\theta^Tg(\theta)\right),-g_1(\theta)\overleftarrow{\varphi'}
    \left(c_\theta^Tg(\theta)\right),\ldots,-g_l(\theta)\overleftarrow{\varphi'}
    \left(c_\theta^Tg(\theta)\right)\right)^T=\underline{0}_l^T.
\end{eqnarray}
and
\begin{equation}\label{condition 2o}
P_0m''(\theta,c_\theta)=P_0\left[-\frac{g_ig_j}
{\varphi''\left(\overleftarrow{\varphi'}
    \left(c_\theta^Tg(\theta)\right)\right)}\right]_{i,j=0,\ldots,l},
\end{equation}
which implies that the matrix $P_0m''(\theta,c_\theta)$ is
symmetric. Under assumption (A.2), by Taylor expansion, there
exists $t_n\in \mathbb{R}^{l+1}$ inside the segment that links
$c_\theta$ and $\widehat{c_\theta}$ with
\begin{equation}\label{Taylor Expansion}
    \begin{array}{ccl}
    0 & = & P_nm'(\theta,\widehat{c_\theta})\\
      & = &
    P_nm'(\theta,c_\theta)+\left(P_n m''(\theta,c_\theta)\right)^T
    \left(\widehat{c_\theta}-c_\theta\right)\\
    & & + \frac{1}{2}\left(\widehat{c_\theta}-c_\theta\right)^T
    P_n m'''(\theta,t_n)
    \left(\widehat{c_\theta}-c_\theta\right),
    \end{array}
\end{equation}
in which, $P_n m'''(\theta,t_n)$ is a $(l+1)-$vector whose entries
are $(l+1)\times (l+1)-$matrices. By (A.2), we have for the
sup-norm of vectors and matrices
\begin{equation*}
    \left\|P_n m'''(\theta,t_n)\right\|:=\left\|\frac{1}{n}
    \sum_{i=1}^n m'''(X_i,\theta,t_n)\right\|\leq
    \frac{1}{n}\sum_{i=1}^{n}|H(X_i)|.
\end{equation*}
By the Law of Large Numbers (LLN), $P_n m'''(\theta,t_n)=O_P(1)$.
So using (A.1), we can write the last term in the right hand side
 of (\ref{Taylor Expansion}) as
$o_P(1)\left(\widehat{c_\theta}-c_\theta\right)$. On the other
hand by (A.3),
$P_nm''(\theta,c_\theta):=\frac{1}{n}\sum_{i=1}^{n}m''(X_i,\theta,c_\theta)$
converges to the matrix $P_0m''(\theta,c_\theta)$. Write
$P_nm''(\theta,c_\theta)$ as $P_0m''(\theta,c_\theta)+o_P(1)$ to
obtain from (\ref{Taylor Expansion})
\begin{equation}\label{equation loi limite}
    -P_nm'(\theta,c_\theta)=\left(P_0m''(\theta,c_\theta)+o_P(1)\right)
    \left(\widehat{c_\theta}-c_\theta\right).
\end{equation}
Under (A.3), by the Central Limit Theorem, we have
$\sqrt{n}P_nm'(\theta,c_\theta)=O_P(1)$, which by (\ref{equation
loi limite}) implies that
$\sqrt{n}\left(\widehat{c_\theta}-c_\theta\right)=O_P(1).$ Hence,
from (\ref{equation loi limite}), we get
\begin{equation}\label{equation qui donne la loi limite}
    \sqrt{n}\left(\widehat{c_\theta}-c_\theta\right)=
    {\left[-P_0m''(\theta,c_\theta)\right]}^{-1}
    \sqrt{n}P_nm'(\theta,c_\theta)+o_P(1).
\end{equation}
Under (A.3), the Central Limit Theorem concludes the proof of part
1. In the case when $P_0$ belongs to $\mathcal{M}_\theta$, then
 $c_\theta^T=(\varphi'(1),\underline{0}^T):=\underline{c}$ and
calculation yields
\begin{equation*}
P_0m'(\theta,\underline{c})m'(\theta,\underline{c})^T=\left(
\begin{array}{cc}
  0 & \underline{0}_l^T \\
  \underline{0}_l & P_0g(\theta)g(\theta)^T \\
\end{array}
\right) ~ \text{ and }~
    -\varphi''(1)P_0m''(\theta,\underline{c})=\left(
\begin{array}{cc}
  1 & \underline{0}_l^T \\
  \underline{0}_l & P_0g(\theta)g(\theta)^T \\
\end{array}
\right).
\end{equation*}
A simple calculation yields (\ref{variance si P est dans Mtheta}).\\

\noindent Proof of part (2). By Taylor expansion, there exists
$\overline{t}_{n}$ inside the segment that links $c_\theta$ and
$~\widehat{c_\theta}$ with
\begin{eqnarray}
\widehat{\phi}_{n}(\mathcal{M}_\theta,P_{0}) & = & P_{n}
m(\theta,\widehat{c_\theta}) \notag \\
& = &
P_{n}m(\theta,c_\theta)+\left(P_{n}m'(\theta,c_\theta)\right)^T
(\widehat{c_\theta}-c_\theta)  \notag \\
& & +\frac{1}{2}{(\widehat{c_\theta}-c_\theta)}^T
\left[P_{n}m''(\theta,c_\theta)\right](\widehat{c_\theta}-c_\theta)  \notag \\
& & + \frac{1}{3!}\sum_{1\leq i,j,k\leq
d}{(\widehat{c_\theta}-c_\theta)}_{i}{(\widehat{c_\theta}-c_\theta)}_{j}\times
\notag \\
& & \times {(\widehat{c_\theta}-c_\theta)}_{k}P_{n}\frac{\partial
^{3}} {\partial t_{i}\partial  t_{j}\partial t_{k}}
m(\theta,\overline{t}_{n}).  \label{taylor1}
\end{eqnarray}
\newline
When $P_0$ belongs to $\mathcal{M}_\theta$, then
$c_\theta^T=\underline{c}$. Hence
$P_{n}m(\theta,c_\theta)=P_{n}m(\theta,\underline{c})=P_n0=0$.
Furthermore, by part (1) in Theorem \ref{lois limite a theta
fixe}, it holds $\sqrt{n}(\widehat{c_\theta}-c_\theta)=O_{p}(1).$
Hence, by (A.1), (A.2) and (A.3), we get
\begin{eqnarray}
\widehat{\phi}_{n}(\mathcal{M}_\theta,P_{0}) & = & \left(
P_{n}m'(\theta,c_\theta)\right)^T(\widehat{c_\theta}-c_\theta)+
\notag \\
& & \frac{1}{2}
{(\widehat{c_\theta}-c_\theta)}^T\left[P_0m''(\theta,c_\theta)\right]
(\widehat{c_\theta}-c_\theta)+o_{P}(1/n), \notag
\end{eqnarray}
which by (\ref{equation qui donne la loi limite}), implies
\begin{eqnarray}
\widehat{\phi}_{n}(\mathcal{M}_\theta,P_{0}) & = & {\left[
P_{n}m'(\theta,c_\theta)\right]}^T{\left[-P_0m''(\theta,c_\theta)\right]}^{-1}
\left[P_{n}m'(\theta,c_\theta)\right] +  \notag \\
& & \frac{1}{2}{\left[P_{n}m'(\theta,c_\theta)\right]}^T{
\left[P_0m''(\theta,c_\theta)\right]}^{-1}
\left[P_{n}m'(\theta,c_\theta)\right]+o_{P}(1/n)  \notag \\
& = & \frac{1}{2}{\left[ P_{n}m'(\theta,c_\theta)\right]}^T{
\left[-P_0m''(\theta,c_\theta)\right]}^{-1} \left[ P_{n}m'(\theta
,c)\right]+o_{P}(1/n). \notag
\end{eqnarray}
This yields to
\begin{equation}
\frac{2n}{\varphi''(1)}\widehat{\phi}_{n}(\mathcal{M}_\theta,P_0)=
{\left[\sqrt{n}P_{n}m'(\theta,c_\theta)\right]}^T
{\left[-\varphi''(1)P_0m''(\theta,c_\theta)\right]
}^{-1}\left[\sqrt{n}P_{n}m'(\theta,c_\theta)\right]+o_{P}(1).
\label{2n Phi n}
\end{equation}
Note that when $P_0$ belongs to $\mathcal{M}_\theta$, then
 $c_\theta^T=\underline{c}$ and
calculation yields
\begin{equation*}
P_0m'(\theta,\underline{c})m'(\theta,\underline{c})^T=\left(
\begin{array}{cc}
  0 & \underline{0}_l^T \\
  \underline{0}_l & P_0g(\theta)g(\theta)^T \\
\end{array}
\right) ~\text{ and }~
    -\varphi''(1)P_0m''(\theta,\underline{c})=\left(%
\begin{array}{cc}
  1 & \underline{0}_l^T \\
  \underline{0}_l & P_0g(\theta)g(\theta)^T \\
\end{array}
\right).
\end{equation*}
 Combining this with (\ref{2n Phi n}), we conclude
the proof of part (2).\\

\noindent Proof of part (3). Since
$\left(\widehat{c_\theta}-c_\theta\right)=O_P(1/\sqrt{n})$ and
$P_nm'(\theta,c_\theta)=P_0m'(\theta,c_\theta)+o_P(1)=0+o_P(1)=o_P(1)$,
then, using (\ref{taylor1}), we obtain
\begin{eqnarray}
\sqrt{n}\left( \widehat{\phi}_{n}(\mathcal{M}_\theta,P_0)-
\phi(\mathcal{M}_\theta,P_0)\right) & = & \sqrt{n}\left(
\widehat{\phi}_{n}(\mathcal{M}_\theta,P_0)-P_0m(\theta,c_\theta)\right)\nonumber\\
 & = &
 \sqrt{n}\left(P_{n}m(\theta,c_\theta)-P_0m(\theta,c_\theta)\right)+o_{P}(1),\nonumber
\end{eqnarray}
and the Central Limit Theorem yields to the conclusion of the
proof of Theorem \ref{lois limite a theta fixe}.

\subsection{Proof of Proposition \ref{consistence de theta n et}}

Define the estimates
\begin{equation*}
\widetilde{\theta}_\phi :=\arg\inf_{\theta\in\Theta}\sup_{t\in
T_\theta} P_nm(\theta,t),
\end{equation*}
\begin{equation*}
\widetilde{\phi}\left(\mathcal{M},P_0\right) :=
\inf_{\theta\in\Theta}\sup_{t\in T_\theta} P_nm(\theta,t)
\end{equation*}
and for all $\theta\in \Theta$,
\begin{equation*}
\widetilde{c_\theta} := \arg\sup_{t\in T_\theta} P_nm(\theta,t).
\end{equation*}
By condition (C.5), for all $n$ sufficiently large, it holds
\begin{equation*}
\widehat{\theta}_\phi=\widetilde{\theta}_\phi ~\text{ and }~
\widehat{\phi}\left(\mathcal{M},P_0\right)=\widetilde{\phi}
\left(\mathcal{M},P_0\right).
\end{equation*}
We prove that $\widetilde{\theta}_\phi$ and
$\widetilde{\phi}\left(\mathcal{M},P_0\right)$ are consistent.
First, we prove the consistency of
$\widetilde{\phi}\left(\mathcal{M},P_0\right)$. We have
\begin{equation*}
    \left|\widetilde{\phi}\left(\mathcal{M},P_0\right)-
    \phi\left(\mathcal{M},P_0\right)\right|=\left|P_nm\left(
    \widetilde{\theta}_\phi,\widetilde{c_{\widetilde{\theta}_\phi}}\right)
    - P_0m(\theta^*,c_{\theta^*})\right|=: |A|.
\end{equation*}
This implies
\begin{equation*}
    P_nm\left(\widetilde{\theta}_\phi,c_{\theta^*}\right)-
    P_0m\left(\widetilde{\theta}_\phi,c_{\theta^*}\right)\leq
    A \leq P_n m\left(\theta^*, \widetilde{c_{\widetilde{\theta}_\phi}}\right)
    -P_0m\left(\theta^*,
    \widetilde{c_{\widetilde{\theta}_\phi}}\right).
\end{equation*}
By condition (C.3), both the RHS and LHS terms in the above
display go to $0$. This implies that $A$ tends to $0$ which
concludes the proof of part (i).\\

\noindent Proof of part (ii). Since for sufficiently large $n$, by
condition (C.5), we have $\widehat{c_\theta}=\widetilde{c_\theta}$
for all $\theta\in\Theta$, the convergence  of
$\sup_{\theta\in\Theta}\left\|\widetilde{c_\theta}-c_\theta\right\|$
to $0$ implies (ii). We prove now that
$\sup_{\theta\in\Theta}\left\|\widetilde{c_\theta}-c_\theta\right\|$
tends to $0$. By the very definition of $\widetilde{c_\theta}$ and
condition (C.3), we have
\begin{eqnarray}
P_nm\left(\theta,\widetilde{c_\theta}\right) & \geq &
P_nm\left(\theta,c_\theta\right)\nonumber\\
& \geq & P_0m\left(\theta,c_\theta\right)-o_P(1),
\end{eqnarray}
where $o_P(1)$ does not depends upon $\theta$ (due to condition
(C.3)). Hence, we have for all $\theta\in\Theta$,
\begin{eqnarray}\label{p 0 m}
P_0m\left(\theta,c_\theta\right)-P_0m\left(\theta,\widetilde{c_\theta}\right)
& \leq &
P_nm\left(\theta,\widetilde{c_\theta}\right)-P_0m\left(\theta,\widetilde{c_\theta}\right)
+o_P(1).
\end{eqnarray}
The term in the RHS of the above display is less than
\begin{equation*}
    \sup_{\theta\in\Theta, t\in T_\theta}\left|P_nm(\theta,t)-
    P_0m(\theta,t)\right|+o_P(1)
\end{equation*}
which by (C.3), tends to $0$. Let $\epsilon > 0$ be such that
$\sup_{\theta\in\Theta}\left\|\widetilde{c_\theta}-c_\theta\right\|>
\epsilon$. There exists some $a_n$ $\in$ $\Theta$ such that
$\left\|\widetilde{c_{a_n}}-c_{a_n}\right\|>\epsilon$. Together
with the strict concavity of the function $t\in T_\theta
\rightarrow P_0m(\theta,t)$ for all $\theta\in\Theta$, there
exists $\eta > 0$ such that
\begin{equation*}
    P_0m\left(a_n,c_{a_n}\right)-P_0m\left(a_n,\widetilde{c_{a_n}}\right)
    > \eta.
\end{equation*}
We then conclude that
\begin{equation*}
    P\left\{\sup_{\theta\in\Theta}\left\|\widetilde{c_\theta}-
    c_\theta\right\|>\epsilon\right\}\leq P\left\{
    P_0m\left(a_n,c_{a_n}\right)-P_0m\left(a_n,\widetilde{c_{a_n}}\right)
    > \eta\right\},
\end{equation*}
and the RHS term tends to $0$ by (\ref{p 0 m}). This concludes the
proof part (ii).\\

\noindent Proof of part (iii). We prove that
$\widetilde{\theta}_\phi$ converges to $\theta^*$. By the very
definition of $\widetilde{\theta}_\phi$, condition (C.4.b) and
part (ii), we obtain
\begin{eqnarray}
P_nm\left(\widetilde{\theta}_\phi,\widetilde{c_{\widetilde{\theta}_\phi}}\right)
& \leq & P_n
m\left(\theta^*,\widetilde{c_{\theta^*}}\right)\nonumber\\ & \leq
&
P_0m\left(\theta^*,\widetilde{c_{\widetilde{\theta}_\phi}}\right)-o_P(1),
\nonumber
\end{eqnarray}
from which
\begin{eqnarray}\label{pour concl}
P_0m\left(\widetilde{\theta}_\phi,\widetilde{c_{\widetilde{\theta}_\phi}}\right)
-P_0m\left(\theta^*,\widetilde{c_{\widetilde{\theta}_\phi}}\right)
& \leq &
P_0m\left(\widetilde{\theta}_\phi,\widetilde{c_{\widetilde{\theta}_\phi}}\right)-
P_nm\left(\widetilde{\theta}_\phi,\widetilde{c_{\widetilde{\theta}_\phi}}\right)
+o_P(1)\nonumber\\
& \leq & \sup_{\left\{\theta\in\Theta, t\in T_\theta\right\}}
\left|P_nm(\theta,t)-P_0m(\theta,t)\right|+o_P(1).
\end{eqnarray}
Further, by part (ii) and condition (C.4.a), for any positive
$\epsilon$, there exists $\eta > 0$ such that
\begin{equation*}
    P\left\{\left\|\widetilde{\theta}_\phi-\theta^*\right\|>\epsilon\right\}
    \leq P\left\{P_0m\left(\widetilde{\theta}_\phi,
    \widetilde{c_{\widetilde{\theta}_\phi}}\right)-P_0
    m\left(\theta^*,
    \widetilde{c_{\widetilde{\theta}_\phi}}\right)>\eta\right\}.
\end{equation*}
The RHS term, under condition (C.3), tends to $0$ by (\ref{pour
concl}). This concludes the proof of Proposition \ref{consistence
de theta n et}.

\subsection{Proof of Theorem \ref{loi limite de thetan et cn}}

Since $P_0\in\mathcal{M}$, then $c_\theta=\underline{c}$. Some
calculus yield
\begin{equation*}
\frac{\partial}{\partial
t}m(\theta_0,\underline{c})={\left[0,-g_1(\theta_0),\ldots,-g_l(\theta_0)\right]}^T=
-{\left[0,{g(\theta_0)}^T\right]}^T,
\end{equation*}
\begin{equation}
\frac{\partial}{\partial\theta}m(\theta,t)=
-\sum_{j=0}^{l}t_i\overleftarrow{\varphi'}
    (t^Tg(\theta))\frac{\partial}{\partial\theta}g_j(\theta), ~~~
\frac{\partial}{\partial
\theta}m(\theta_0,\underline{c})=\underline{0}_d,
\end{equation}
\begin{equation*}
\frac{\partial^2}{\partial\theta\partial
t}m(\theta_0,\underline{c})=\left[\underline{0}_d,-\frac{\partial}{\partial\theta}
g_1(\theta_0),\ldots,-\frac{\partial}{\partial\theta}g_l(\theta_0)\right]
=-\left[\underline{0}_d,\frac{\partial}{\partial\theta}g(\theta)\right]
\end{equation*}
\begin{equation*}
\frac{\partial^2}{\partial
t\partial\theta}m(\theta_0,\underline{c})=\left[\frac{\partial^2}{\partial\theta\partial
t}m(\theta_0,\underline{c})\right]^T,~
   \frac{\partial^2}{\partial\theta^2}m(\theta_0,\underline{c})=
   \left[\underline{0}_d,\ldots,\underline{0}_d\right],
\end{equation*}
and
\begin{equation*}
    \frac{\partial^2}{\partial t^2}m(\theta_0,\underline{c})=\frac{1}{\varphi''(1)}
    {\left[-g_i(\theta_0)g_j(\theta_0)\right]}_{i,j=0,1,\ldots,l}:=
    \frac{-1}{\varphi''(1)}\left(\overline{g}(\theta_0)
    \overline{g}(\theta_0)^T\right).
\end{equation*}
Integrating w.r.t. $P_0$, we obtain
\begin{equation}
P_0\frac{\partial}{\partial
t}m(\theta_0,\underline{c})=\underline{0}_l,~~
P_0\frac{\partial}{\partial
\theta}m(\theta_0,\underline{c})=\underline{0}_d,
~P_0\frac{\partial^2}{\partial\theta^2}m(\theta_0,\underline{c})=
\left[\underline{0}_d,\ldots,\underline{0}_d\right],
\end{equation}
\begin{equation}
    P_0\frac{\partial^2}{\partial\theta\partial
    t}m(\theta_0,\underline{c})=
    -\left[\underline{0}_d,P_0\frac{\partial}{\partial\theta}g(\theta_0)\right],
\end{equation}
\begin{equation}
P_0\frac{\partial^2}{\partial
t\partial\theta}m(\theta_0,\underline{c})= \left[
P_0\frac{\partial^2}{\partial\theta\partial
    t}m(\theta_0,\underline{c})\right]^T=-
    \left[\underline{0}_d,P_0\frac{\partial}{\partial\theta}g(\theta_0)\right]^T,
\end{equation}
and
\begin{eqnarray}
  P_0    \frac{\partial^2}{\partial
  t^2}m(\theta_0,\underline{c}) & = &
  \frac{-1}{\varphi''(1)}\left[P_0g_i(\theta_0)g_j(\theta_0)
  \right]_{i,j=0,1,\ldots,l}.\nonumber\\
  ~ & = & \frac{-1}{\varphi''(1)}\left(%
\begin{array}{cc}
  1 & \underline{0}_l^T \\
  \underline{0}_l & P_0g(\theta_0)g(\theta_0)^T \\
\end{array}%
\right)
\end{eqnarray}
By the very definition of $\widehat{\theta}_\phi$ and
$\widehat{c_{\widehat{\theta}_\phi}}$, they both obey
\begin{equation*}
\left\{
\begin{array}{ccc}
P_n\frac{\partial}{\partial t}m\left(\theta,t\right)
& = & 0\\
P_n\frac{\partial}{\partial\theta} m\left(\theta,t(\theta)\right)
& = & 0,
\end{array}
\right.
\end{equation*}
i.e.,
\begin{equation*}
\left\{
\begin{array}{ccc}
P_n\frac{\partial}{\partial
t}m\left(\widehat{\theta}_\phi,\widehat{c_{\widehat{\theta}_\phi}}\right)
& = & 0\\
P_n\frac{\partial}{\partial\theta}
m\left(\widehat{\theta}_\phi,\widehat{c_{\widehat{\theta}_\phi}}\right)+P_n
\frac{\partial}{\partial t}
m\left(\widehat{\theta}_\phi,\widehat{c_{\widehat{\theta}_\phi}}\right)
\frac{\partial}{\partial
\theta}\widehat{c_{\widehat{\theta}_\phi}} & = & 0.
\end{array}
\right.
\end{equation*}
The second term in the left hand side of the second equation is
equal to $0$, due to the first equation. Hence
$\widehat{c_{\widehat{\theta}_\phi}}$ and $\widehat{\theta}_\phi$
are solutions of the somehow simpler system
\begin{equation*}
\left\{
\begin{array}{ccc}
P_n\frac{\partial}{\partial
t}m\left(\widehat{\theta}_\phi,\widehat{c_{\widehat{\theta}_\phi}}\right)
& = & 0~~~~~~~~~~~~~~~~(E1)\\
P_n\frac{\partial}{\partial\theta}
m\left(\widehat{\theta}_\phi,\widehat{c_{\widehat{\theta}_\phi}}\right)
& = & 0~~~~~~~~~~~~~~(E2).
\end{array}
\right.
\end{equation*}
Use a Taylor expansion in (E1); there exists
$\left(\widetilde{\theta}_n,\widetilde{t}_n\right)$ inside the
segment that links
$(\widehat{\theta}_\phi,\widehat{c_{\widehat{\theta}_\phi}})$ and
$(\theta_0,\underline{c})$ such that
\begin{eqnarray}\label{Taylor 1 E1}
  0 & = & P_n\frac{\partial}{\partial t}m\left(\theta_0,\underline{c}\right)+
  \left[\left(P_n\frac{\partial^2}{\partial t^2}m(\theta_0,\underline{c})\right)^T,
  \left(P_n\frac{\partial^2}{\partial \theta\partial t}
  m(\theta_0,\underline{c})\right)^T\right]a_n\nonumber\\
  & & +\frac{1}{2}a_n^TA_na_n,
\end{eqnarray}
with
\begin{equation}\label{a n }
    a_n:={\left({\left(\widehat{c_{\widehat{\theta}_\phi}}-\underline{c}\right)}^T,
{\left(\widehat{\theta}_\phi-\theta_0\right)}^T\right)}^T
\end{equation}
and
\begin{equation}
    A_n:=\left(%
\begin{array}{cc}
  P_n\frac{\partial ^3}{\partial t ^3}m(\widetilde{\theta},\widetilde{c}_n) &
  P_n \frac{\partial ^3}{\partial t\partial \theta\partial t}
  m(\widetilde{\theta},\widetilde{c}_n)\\
  P_n\frac{\partial ^3}{\partial \theta\partial t ^2}m(\widetilde{\theta},\widetilde{c}_n) &
  P_n\frac{\partial ^3}{\partial\theta^2\partial t}m(\widetilde{\theta},\widetilde{c}_n) \\
\end{array}
\right).
\end{equation}
By (A.5), the LLN implies that $A_n=O_P(1)$. So using (A.4), we
can write the last term in right hand side of (\ref{Taylor 1 E1})
as $o_P(1)a_n$. On the other hand by (A.6), we can write also\\
$\left[\left(P_n\frac{\partial^2}{\partial
t^2}m(\theta_0,\underline{c})\right)^T,
  \left(P_n\frac{\partial^2}{\partial\theta\partial t}
  m(\theta_0,\underline{c})\right)^T\right]$ as
  $\left[P_0\frac{\partial^2}{\partial t^2}m(\theta_0,\underline{c}),
  \left(P_0\frac{\partial^2}{\partial\theta\partial t}
  m(\theta_0,\underline{c})\right)^T\right]+o_P(1)$ to obtain
  from (\ref{Taylor 1 E1})
  \begin{equation}\label{Equation E1 donne}
    -P_n\frac{\partial}{\partial
    t}m(\theta_0,\underline{c})=
    \left[P_0\frac{\partial^2}{\partial t^2}m(\theta_0,\underline{c})+o_P(1),
    \left(P_0\frac{\partial^2}{\partial\theta\partial
    t}m(\theta_0,\underline{c})\right)^T
    +o_P(1)\right]a_n.
\end{equation}
In the same way, using a Taylor expansion in (E2), there exists
$(\overline{\theta}_n,\overline{t}_n)$ inside the segment that
links
$\left(\widehat{\theta}_\phi,\widehat{c_{\widehat{\theta}_\phi}}\right)$
and $(\theta_0,\underline{c})$ such that
\begin{eqnarray}\label{Taylor 1 E2}
0 & = & P_{n}\frac{\partial }{\partial \theta }m(\theta
_{0},\underline{c})+\left[ \left( P_{n}\frac{\partial
^{2}}{\partial t
\partial \theta }m(\theta_{0},\underline{c})\right)^T,\left(
P_{n}\frac{\partial ^{2}}{\partial
\theta^2}m(\theta_{0},\underline{c})\right)^T\right] a_{n}  \notag \\
& & +\frac{1}{2}a_{n}^{t}B_{n}a_{n},
\end{eqnarray}
with
\begin{equation*}
B_n:=\left[
\begin{array}{cc}
P_{n}\frac{\partial ^{3}}{\partial t^2\partial
\theta}m(\overline{\theta}_{n},\overline{t}_{n}) &
P_{n}\frac{\partial^3}{\partial t \partial \theta^2
}m(\overline{\theta}_n,\overline{t}_n)
\\
P_{n}\frac{\partial^3}{\partial \theta \partial t \partial \theta
}m(\overline{\theta}_n,\overline{t}_n) & P_{n}\frac{\partial
^{3}}{\partial \theta^3}m(\overline{\theta}_n,\overline{t}_{n})
\end{array}
\right].
\end{equation*}
As in (\ref{Equation E1 donne}), we obtain
  \begin{equation}\label{Equation E2 donne}
    -P_n\frac{\partial}{\partial
    \theta}m(\theta_0,\underline{c})=
    \left[\left(P_0\frac{\partial^2}{\partial t\partial\theta}m(\theta_0,\underline{c})
    \right)^T+o_P(1),
    P_0\frac{\partial^2}{\partial
  \theta^2}m(\theta_0,\underline{c})+o_P(1)\right]a_n.
\end{equation}
From (\ref{Equation E1 donne}) and (\ref{Equation E2 donne}), we
get
\begin{eqnarray}\label{racine de n a n 1}
\sqrt{n}a_n  & = &  \sqrt{n}\left(
\begin{array}{cc}
  P_0\frac{\partial ^2}{\partial t^2}m(\theta_0,c(\theta_0)) & \left(P_0\frac{\partial^2}
  {\partial\theta\partial t} m(\theta_0,\underline{c})\right)^T\\
  \left(P_0\frac{\partial^2}{\partial t\partial \theta} m(\theta_0,c(\theta_0))\right)^T &
  P_0\frac{\partial^2}{\partial \theta^2}m(\theta_0,c(\theta_0)) \\
\end{array}
\right)^{-1}\times\nonumber\\
& & \times\left(
\begin{array}{c}
  -P_n\frac{\partial}{\partial t}m(\theta_0,\underline{c}) \\
  -P_n\frac{\partial }{\partial \theta} m(\theta_0,\underline{c})\\
\end{array}
\right)+o_P(1).
\end{eqnarray}
Denote $S$ the $(l+1+d)\times(l+1+d)-$matrix defined by
\begin{equation}\label{S}
    S:=\left(
\begin{array}{cc}
  S_{11} & S_{12} \\
  S_{21} & S_{22} \\
\end{array}
\right):=\left(
\begin{array}{cc}
  P_0\frac{\partial ^2}{\partial t^2}m(\theta_0,c(\theta_0)) & \left(P_0\frac{\partial^2}
  {\partial\theta\partial t} m(\theta_0,\underline{c})\right)^T\\
  \left(P_0\frac{\partial^2}{\partial t\partial \theta} m(\theta_0,c(\theta_0))\right)^T &
  P_0\frac{\partial^2}{\partial \theta^2}m(\theta_0,c(\theta_0)) \\
\end{array}
\right).
\end{equation}
We have
\begin{equation}
S_{11}=\frac{-1}{\varphi''(1)}\left(
\begin{array}{cc}
  1 & \underline{0}_l^T \\
  \underline{0}_l & P_0g(\theta_0)g(\theta_0)^T \\
\end{array}
\right)
\end{equation}
\begin{equation}
S_{12}=-\left[\underline{0}_d,
P_0\frac{\partial}{\partial\theta}g(\theta_0)\right]^T,~
S_{21}=-\left[\underline{0}_d,
P_0\frac{\partial}{\partial\theta}g(\theta_0)\right]\text{ and }
\end{equation}
\begin{equation}
S_{22}=P_0\frac{\partial^2}{\partial\theta^2}m(\theta_0,\underline{c})
=\left[\underline{0}_d,\ldots,\underline{0}_d\right].
\end{equation}
The inverse matrix $S^{-1}$ of the matrix $S$ writes
\begin{equation}\label{S - 1}
S^{-1}=\left(
\begin{array}{cc}
  S_{11}^{-1}+S_{11}^{-1}S_{12}S_{22.1}^{-1}S_{21}S_{11}^{-1} &
  -S_{11}^{-1}S_{12}S_{22.1}^{-1} \\
  -S_{22.1}^{-1}S_{21}S_{11}^{-1} & S_{22.1}^{-1} \\
\end{array}
\right),
\end{equation}
where
\begin{eqnarray}
S_{22.1} & = & -S_{21}S_{11}^{-1}S_{12}\nonumber\\
         & = & \left[\underline{0}_d,
P_0\frac{\partial}{\partial\theta}g(\theta_0)\right]\left[\varphi''(1)\right]
\left[
\begin{array}{cc}
  1 & \underline{0}_l^T \\
  \underline{0}_l & \left[P_0g(\theta_0)g(\theta_0)^T\right]^{-1} \\
\end{array}
\right] \left[\underline{0}_d,
P_0\frac{\partial}{\partial\theta}g(\theta_0)\right]^T\nonumber\\
& = & \varphi''(1)
    \left[P_0\frac{\partial}{\partial\theta}g(\theta_0)\right]
    \left[P_0g(\theta_0)g(\theta_0)^T\right]^{-1}
    \left[P_0\frac{\partial}{\partial\theta}g(\theta_0)\right]^T.
\end{eqnarray}

From (\ref{racine de n a n 1}), using (\ref{S}) and (\ref{S - 1}),
we can write
\begin{eqnarray}\label{racine de n a n 3}
\sqrt{n}\left(
\begin{array}{c}
  \widehat{c_{\widehat{\theta}_\phi}}-\underline{c} \\
  \widehat{\theta}_\phi-\theta_0 \\
\end{array}
\right) & = & \left(
\begin{array}{cc}
  S_{11}^{-1}+S_{11}^{-1}S_{12}S_{22.1}^{-1}S_{21}S_{11}^{-1} &
  -S_{11}^{-1}S_{12}S_{22.1}^{-1} \\
  -S_{22.1}^{-1}S_{21}S_{11}^{-1} & S_{22.1}^{-1} \\
\end{array}
\right)\times\nonumber\\
& & \times\sqrt{n}\left(
\begin{array}{c}
  -P_n\frac{\partial}{\partial t}m(\theta_0,\underline{c}) \\
  \underline{0}_d \\
\end{array}
\right)+o_P(1).
\end{eqnarray}
Note that
\begin{equation}\label{racine de n - Pn 0d}
\sqrt{n}\left(
\begin{array}{c}
  -P_n\frac{\partial}{\partial t}m(\theta_0,\underline{c}) \\
  \underline{0}_d \\
\end{array}
\right),
\end{equation}
under assumption (A.6), by the Central Limit Theorem, converges in
distribution to a centered multivariate normal variable with
covariance matrix
\begin{equation}\label{M}
    M=\left(
\begin{array}{cc}
  M_{11} & M_{12} \\
  M_{21} & M_{22} \\
\end{array}
\right)
\end{equation}
where
\begin{equation}\label{M i j}
    M_{11}=\left(
\begin{array}{cc}
  0   & \underline{0}_l^T \\
  \underline{0}_l & P_0g(\theta_0)g(\theta_0)^T \\
\end{array}
\right),~ M_{12}= \left(
\begin{array}{c}
  \underline{0}_d^T \\
  \vdots \\
  \underline{0}_d^T \\
\end{array}
\right),~
    M_{21}=\left(
\begin{array}{cc}
  0 & \underline{0}_l^T \\
  \vdots & \vdots \\
  0 &  \underline{0}_l^T\\
\end{array}
\right) ~\text{ and }~
    M_{22}=\left(
\begin{array}{c}
  \underline{0}_d^T \\
  \vdots \\
  \underline{0}_d^T \\
\end{array}
\right).
\end{equation}
Hence, from (\ref{racine de n - Pn 0d}), we deduce that
\begin{equation}\label{racine de n c n - c theta n - theta 0}
\sqrt{n}\left(
\begin{array}{c}
  \widehat{c_{\widehat{\theta}_\phi}}-\underline{c} \\
  \widehat{\theta}_\phi-\theta_0 \\
\end{array}
\right)
\end{equation}
converges in distribution to a centered multivariate normal
variable with covariance matrix
\begin{equation}\label{CM nonparam}
    C=S^{-1}M {\left[S^{-1}\right]}^T:=\left(
\begin{array}{cc}
  C_{11} & C_{12} \\
  C_{21} & C_{22} \\
\end{array}
\right),
\end{equation}
and using (\ref{M i j}) and some algebra, we get
\begin{eqnarray}\label{C 11}
    C_{11} & = & {\varphi''(1)}^2\left[
\begin{array}{cc}
  0 & \underline{0}_l^T \\
  \underline{0}_l & \left[P_0g(\theta_0)g(\theta_0)^T\right]^{-1} \\
\end{array}
\right]-{\varphi''(1)}^2\left[
\begin{array}{cc}
  0 & \underline{0}_l^T \\
  \underline{0}_l & \left[P_0g(\theta_0)g(\theta_0)^T\right]^{-1} \\
\end{array}
\right]\times\nonumber\\
& &
\times\left[\underline{0},P_0\frac{\partial}{\partial\theta}g(\theta_0)\right]^T
\left[C_{22}\right]
\left[\underline{0},P_0\frac{\partial}{\partial\theta}g(\theta_0)\right]\left[
\begin{array}{cc}
  0 & \underline{0}_l^T \\
  \underline{0}_l & \left[P_0g(\theta_0)g(\theta_0)^T\right]^{-1} \\
\end{array}
\right],
\end{eqnarray}
\begin{equation}\label{C 12 et C 21}
    C_{12}=\left[\underline{0}_l,\ldots,\underline{0}_l\right],~
    C_{21}=\left[\underline{0}_d,\ldots,\underline{0}_d\right]
\end{equation}
and
\begin{equation}\label{C 22}
    C_{22}={\left\{\left[P_0\frac{\partial}{\partial\theta}g(\theta_0)\right]
  \left[P_0\left(g(\theta_0)g(\theta_0)^T\right)\right]^{-1}\left[P_0\frac{\partial
  }{\partial\theta}g(\theta_0)\right]^T\right\}}^{-1}.
\end{equation}
From (\ref{racine de n c n - c theta n - theta 0}), we deduce that
$C_{11}$ and $C_{22}$ are respectively the limit covariance matrix
of
$\sqrt{n}\left(\widehat{c_{\widehat{\theta}_\phi}}-\underline{c}\right)$
and $\sqrt{n}\left(\widehat{\theta}_\phi-\theta_0\right)$, i.e.,
$U=C_{11}$ and $V=C_{22}$. (\ref{C 12 et C 21}) implies that
$\sqrt{n}\left(\widehat{c_{\widehat{\theta}_\phi}}-\underline{c}\right)$
and $\sqrt{n}\left(\widehat{\theta}_\phi-\theta_0\right)$ are
asymptotically uncorrelated. This concludes the Proof of Theorem
\ref{loi limite de thetan et cn}.

\subsection{Proof of Theorem \ref{loi limite du vect thetan  cn}}

Under assumptions (A.4-6), as in the proof of Theorem \ref{loi
limite de thetan et cn}, we obtain
\begin{equation*}
\sqrt{n}\left(%
\begin{array}{c}
  \widehat{c_{\widehat{\theta}_\phi}}-c_{\theta^*} \\
  \widehat{\theta}_\phi -\theta^*\\
\end{array}%
\right)=\sqrt{n}S^{-1}\left(%
\begin{array}{c}
  -P_n\frac{\partial}{\partial t}m(\theta^*,c_{\theta^*}) \\
  -P_n\frac{\partial}{\partial \theta}m(\theta^*,c_{\theta^*})\\
\end{array}%
\right)+o_P(1),
\end{equation*}
and the CLT concludes the proof.

\subsection{Proof of Theorem \ref{loi limite de F chapeau n}}

Using Taylor expansion at $(\underline{c},\theta_0)$, we get
\begin{eqnarray}\label{Taylor F chapeau n}
 \widehat{F}_n(x) & := & \sum_{i=1}^n\widehat{Q^*_{\widehat{\theta}_\phi}}
 \mathds{1}_{(-\infty,x]}(X_i):=\frac{1}{n}\sum_{i=1}^n\overleftarrow{\varphi'}
 \left(\widehat{c_{\widehat{\theta}_\phi}}^T\overline{g}(X_i,\widehat{\theta}_\phi)
 \right)\mathds{1}_{(-\infty,x]}(X_i)\nonumber\\
 & = & F_n(x)+
   \frac{1}{n}\left[\sum_{i=1}^n \overline{g}(X_i,\theta_0)
   \mathds{1}_{(-\infty,x]}(X_i)\right]^T
   \frac{1}{\varphi''(1)}\left(\widehat{c_{\widehat{\theta}_\phi}}-\underline{c}\right)
   +o_P(\delta_n),\nonumber\\
   & &
\end{eqnarray}
where
$\delta_n:=\left\|\widehat{c_{\widehat{\theta}_\phi}}-\underline{c}\right\|+
\left\|\widehat{\theta}_\phi-\theta_0\right\|$, which by Theorem
\ref{loi limite de thetan et cn}, is equal to $O_P(1/\sqrt{n})$.
Hence, (\ref{Taylor F chapeau n}) yields
\begin{eqnarray}\label{Taylor 2 F chapeau n}
  \sqrt{n}\left(\widehat{F}_n(x)-F(x)\right) & = &
  \sqrt{n}\left(F_n(x)-F(x)\right)+\nonumber\\
    & & +\frac{1}{\varphi''(1)}
  \left[P_0\left(\overline{g}(\theta_0)\mathds{1}_{(-\infty,x]}\right)\right]^T
  \sqrt{n}\left(\widehat{c_{\widehat{\theta}_\phi}}-\underline{c}\right)
   +o_P(1).\nonumber\\
   & &
\end{eqnarray}
On the other hand, from (\ref{racine de n a n 3}), we get
\begin{eqnarray}\label{racine de n cn -c}
\sqrt{n}\left(\widehat{c_{\widehat{\theta}_\phi}}-\underline{c}\right)
 & = & H\sqrt{n}\left(-P_n\frac{\partial}{\partial\theta}
 m(\theta_0,\underline{c})\right)+o_P(1)
\end{eqnarray}
with
\begin{equation}\label{H}
    H=S_{11}^{-1}+S_{11}^{-1}S_{12}S_{22.1}^{-1}S_{21}S_{11}^{-1}.
\end{equation}
We will use $f(.)$ to denote the function
$\mathds{1}_{(-\infty,x]}(.)-F(x)$, for all $x\in\mathbb{R}$.
Substituting (\ref{racine de n cn -c}) in (\ref{Taylor 2 F chapeau
n}), we get
\begin{eqnarray}\label{Taylor 3 F chapeau n}
  \sqrt{n}\left(\widehat{F}_n(x)-F(x)\right) & = &
  \sqrt{n}P_nf+\frac{1}{\varphi''(1)}
  \left[P_0\left(\overline{g}(\theta_0)\mathds{1}_{(-\infty,x]}\right)\right]^T
  H\times\nonumber\\
    & & \times\sqrt{n}\left(-P_n\frac{\partial}{\partial\theta}
 m(\theta_0,\underline{c})\right)
   +o_P(1).
\end{eqnarray}
By the Multivariate Central Limit Theorem, the vector\\ $\sqrt{n}
\left(P_nf,-\left[P_n\frac{\partial}{\partial\theta}
 m(\theta_0,\underline{c})\right]^T\right)^T$
converges in distribution to a centered multivariate normal
variable which implies that
$\sqrt{n}\left(\widehat{F}_n(x)-F(x)\right)$ is asymptotically
centered normal variable. We calculate now  its limit variance,
noted  $W(x)$.
\begin{eqnarray}
 W(x) & = &
 F(x)(1-F(x))+\frac{1}{\varphi''(1)^2}\left[P_0\left(\overline{g}(\theta_0)
 \mathds{1}_{(-\infty,x]}\right)\right]^T U\left[P_0\left(\overline{g}(\theta_0)
 \mathds{1}_{(-\infty,x]}\right)\right]+\nonumber\\
  & & + 2 \frac{1}{\varphi''(1)}\left[P_0\left(\overline{g}(\theta_0)
 \mathds{1}_{(-\infty,x]}\right)\right]^T H
 \left[P_0\left(-\frac{\partial}{\partial t}m(\theta_0,\underline{c})
 \mathds{1}_{(-\infty,x]}\right)\right].
\end{eqnarray}
Use the explicit forms of $\frac{\partial}{\partial
t}m(\theta_0,\underline{c})$, the matrices $U$ and $V$ and some
algebra to obtain
\begin{equation*}
    W(x)=F(x)\left(1-F(x)\right)-
    {\left[P_0\left(g(\theta_0)\mathds{1}_{(-\infty,x]}\right)\right]}^T
    \Gamma
    \left[P_0\left(g(\theta_0)\mathds{1}_{(-\infty,x]}\right)\right],
\end{equation*}
with
\begin{eqnarray}
 \Gamma & = & {\left[P_0g(\theta_0)g(\theta_0)^T\right]}^{-1}
 - {\left[P_0g(\theta_0)g(\theta_0)^T\right]}^{-1}{\left[P_0
 \frac{\partial}{\partial\theta}g(\theta_0)\right]}^T V\times\nonumber\\
     & & \times \left[P_0\frac{\partial}{\partial\theta}g(\theta_0)\right]
 {\left[P_0g(\theta_0)g(\theta_0)^T\right]}^{-1}.\nonumber
\end{eqnarray}
This concludes the proof of Theorem \ref{loi limite de F chapeau
n}.


\begin{thebibliography}{}

\bibitem[Baggerly(1998)]{Baggerly1998}
Baggerly, K.~A. (1998).
\newblock Empirical likelihood as a goodness-of-fit measure.
\newblock {\em Biometrika}, {\bf 85}(3), 535--547.

\bibitem[Basu and Lindsay(1994)]{Basu-Lindsay1994}
Basu, A. and Lindsay, B.~G. (1994).
\newblock Minimum disparity estimation for continuous models: efficiency,
  distributions and robustness.
\newblock {\em Ann. Inst. Statist. Math.}, {\bf 46}(4), 683--705.

\bibitem[Bertail(2004)]{Bertail2003}
Bertail, P. (2004).
\newblock Empirical likelihood in nonparametric and semiparametric models.
\newblock In {\em Parametric and semiparametric models with applications to
  reliability, survival analysis, and quality of life}, Stat. Ind. Technol.,
  pages 291--306. Birkh\"auser Boston, Boston, MA.



\bibitem[Bertail(2006)]{Bertail2006}
Bertail, P. (2006).
\newblock Empirical likelihood in some semiparametric models
\newblock {\em Bernoulli}, {\bf 12}(2), 299--331.




\bibitem[Bickel {\em et~al.}(1991)]{BickelRiovWellner1991}
Bickel, P.~J., Ritov, Y., and Wellner, J.~A. (1991).
\newblock Efficient estimation of linear functionals of a probability measure
  {$P$} with known marginal distributions.
\newblock {\em Ann. Statist.}, {\bf 19}(3), 1316--1346.

\bibitem[Borwein and Lewis(1991)]{BorweinLewis1991}
Borwein, J.~M. and Lewis, A.~S. (1991).
\newblock Duality relationships for entropy-like minimization problems.
\newblock {\em SIAM J. Control Optim.}, {\bf 29}(2), 325--338.

\bibitem[Borwein and Lewis(1993)]{BorweinLewis1993}
Borwein, J.~M. and Lewis, A.~S. (1993).
\newblock Partially-finite programming in {$L\sb 1$} and the existence of
  maximum entropy estimates.
\newblock {\em SIAM J. Optim.}, {\bf 3}(2), 248--267.

\bibitem[Broniatowski and Keziou(2006)]{Bronia_Kez2006_STUDIA}
Broniatowski, M. and Keziou, A. (2006).
\newblock Minimization of {$\phi$}-divergences on sets of signed measures.
\newblock {\em Studia Sci. Math. Hungar.}, {\bf 43}(4), 403--442.

\bibitem[Broniatowski and Keziou(2009)]{Broniatowski-Keziou2008}
Broniatowski, M. and Keziou, A. (2009).
\newblock Parametric estimation and testing through divergences.
\newblock {\em \textit{ Journal of Multivariate Analysis}}, {\bf
43}, 16-36.

\bibitem[Corcoran(1998)]{Corcoran1998}
Corcoran, S. (1998).
\newblock Bertlett adjustement of empirical discrepancy statistics.
\newblock {\em Biometrika}, {\bf 85}, 967--972.

\bibitem[Cressie and Read(1984)]{Cressie-Read1984}
Cressie, N. and Read, T. R.~C. (1984).
\newblock Multinomial goodness-of-fit tests.
\newblock {\em J. Roy. Statist. Soc. Ser. B}, {\bf 46}(3), 440--464.

\bibitem[Csisz{\'a}r(1963)]{Csiszar1963}
Csisz{\'a}r, I. (1963).
\newblock Eine informationstheoretische {U}ngleichung und ihre {A}nwendung auf
  den {B}eweis der {E}rgodizit\"at von {M}arkoffschen {K}etten.
\newblock {\em Magyar Tud. Akad. Mat. Kutat\'o Int. K\"ozl.}, {\bf 8}, 85--108.

\bibitem[Csisz{\'a}r(1967)]{Csiszar1967}
Csisz{\'a}r, I. (1967).
\newblock On topology properties of {$f$}-divergences.
\newblock {\em Studia Sci. Math. Hungar.}, {\bf 2}, 329--339.

\bibitem[Csisz{\'a}r(1975)]{Csiszar1975}
Csisz{\'a}r, I. (1975).
\newblock {$I$}-divergence geometry of probability distributions and
  minimization problems.
\newblock {\em Ann. Probability}, {\bf 3}, 146--158.

\bibitem[Haberman(1984)]{Haberman1984}
Haberman, S.~J. (1984).
\newblock Adjustment by minimum discriminant information.
\newblock {\em Ann. Statist.}, {\bf 12}(3), 971--988.

\bibitem[Hansen(1982)]{Hansen1982}
Hansen, L.~P. (1982).
\newblock Large sample properties of generalized method of moments estimators.
\newblock {\em Econometrica}, {\bf 50}(4), 1029--1054.

\bibitem[Imbens(1997)]{Imbens1997}
Imbens, G.~W. (1997).
\newblock One-step estimators for over-identified generalized method of moments
  models.
\newblock {\em Rev. Econom. Stud.}, {\bf 64}(3), 359--383.

\bibitem[Jim{\'e}nez and Shao(2001)]{JimenezShao2001}
Jim{\'e}nez, R. and Shao, Y. (2001).
\newblock On robustness and efficiency of minimum divergence estimators.
\newblock {\em Test}, {\bf 10}(2), 241--248.

\bibitem[Keziou(2003)]{Keziou2003}
Keziou, A. (2003).
\newblock Dual representation of {$\phi$}-divergences and applications.
\newblock {\em C. R. Math. Acad. Sci. Paris}, {\bf 336}(10), 857--862.

\bibitem[L{\'e}onard(2001a)]{Leonard2001a}
L{\'e}onard, C. (2001a).
\newblock Convex conjugates of integral functionals.
\newblock {\em Acta Math. Hungar.}, {\bf 93}(4), 253--280.

\bibitem[L{\'e}onard(2001b)]{Leonard2001b}
L{\'e}onard, C. (2001b).
\newblock Minimizers of energy functionals.
\newblock {\em Acta Math. Hungar.}, {\bf 93}(4), 281--325.

\bibitem[Liese and Vajda(1987)]{Liese-Vajda1987}
Liese, F. and Vajda, I. (1987).
\newblock {\em Convex statistical distances}, volume~95.
\newblock BSB B. G. Teubner Verlagsgesellschaft, Leipzig.

\bibitem[Lindsay(1994)]{Lindsay1994}
Lindsay, B.~G. (1994).
\newblock Efficiency versus robustness: the case for minimum {H}ellinger
  distance and related methods.
\newblock {\em Ann. Statist.}, {\bf 22}(2), 1081--1114.

\bibitem[Newey and Smith(2004)]{NeweySmith2003}
Newey, W.~K. and Smith, R.~J. (2004).
\newblock Higher order properties of \text{GMM} and generalized empirical
  likelihood estimators.
\newblock {\em Econometrica}.

\bibitem[Neyman(1937)]{Neyman1937}
Neyman, J. (1937).
\newblock Outline of a theory of statistical estimation based on the classical
  theory of probability.
\newblock {\em Phil. Trans. Roy. Soc. Ser.}, {\bf A}(236), 333--380.

\bibitem[Nikitin(1995)]{Nikitin1995}
Nikitin, Y. (1995).
\newblock {\em Asymptotic efficiency of nonparametric tests}.
\newblock Cambridge University Press, Cambridge.

\bibitem[Owen(1990)]{Owen1990}
Owen, A. (1990).
\newblock Empirical likelihood ratio confidence regions.
\newblock {\em Ann. Statist.}, {\bf 18}(1), 90--120.

\bibitem[Owen(1988)]{Owen1988}
Owen, A.~B. (1988).
\newblock Empirical likelihood ratio confidence intervals for a single
  functional.
\newblock {\em Biometrika}, {\bf 75}(2), 237--249.

\bibitem[Owen(2001)]{Owen2001}
Owen, A.~B. (2001).
\newblock {\em Empirical Likelihood}.
\newblock Chapman and Hall, New York.

\bibitem[Qin and Lawless(1994)]{Qin-Lawless1994}
Qin, J. and Lawless, J. (1994).
\newblock Empirical likelihood and general estimating equations.
\newblock {\em Ann. Statist.}, {\bf 22}(1), 300--325.

\bibitem[Rao(1961)]{Rao1961}
Rao, C.~R. (1961).
\newblock Asymptotic efficiency and limiting information.
\newblock In {\em Proc. 4th Berkeley Sympos. Math. Statist. and Prob., Vol. I},
  pages 531--545. Univ. California Press, Berkeley, Calif.

\bibitem[R{\"u}schendorf(1984)]{Ruschendorf1984}
R{\"u}schendorf, L. (1984).
\newblock On the minimum discrimination information theorem.
\newblock {\em Statist. Decisions}, (suppl. 1), 263--283.
\newblock Recent results in estimation theory and related topics.

\bibitem[Schennach(2007)]{Schennach2007}
Schennach, S.~M. (2007).
\newblock Point estimation with exponentially tilted empirical likelihood.
\newblock {\em Ann. Statist.}, {\bf 35}(2), 634--672.

\bibitem[Sen and Singer(1993)]{Sen-Singer1993}
Sen, P.~K. and Singer, J.~M. (1993).
\newblock {\em Large sample methods in statistics}.
\newblock Chapman \& Hall, New York.

\bibitem[Sheehy(1987)]{Sheehy1987}
Sheehy, A. (1987).
\newblock \text{Kullback-Leibler} constrained estimation of probability
  measures.
\newblock {\em Report, Dept. Statistics, Stanford Univ.}

\bibitem[Takagi(1998)]{Takagi1998}
Takagi, Y. (1998).
\newblock A new criterion of confidence set estimation: improvement of the
  {N}eyman shortness.
\newblock {\em J. Statist. Plann. Inference}, {\bf 69}(2), 329--338.

\bibitem[van~der Vaart(1998)]{vanderVaart1998}
van~der Vaart, A.~W. (1998).
\newblock {\em Asymptotic statistics}.
\newblock Cambridge Series in Statistical and Probabilistic Mathematics.
  Cambridge University Press, Cambridge.

\bibitem[Wilks(1938)]{Wilks1938}
Wilks, S. (1938).
\newblock The large-sample distribution of the likelihood ratio for testing
  composite hypotheses.
\newblock {\em Ann. Math. Statist.}, {\bf 9}, 60--62.

\end{thebibliography}

\end{document}